\documentclass{elsart3-1}
\usepackage{amssymb}
\usepackage[english,francais]{babel}
\usepackage{graphicx}
\usepackage{amsmath}
\usepackage{subfigure}
\usepackage{esint}
\usepackage{wasysym}
\usepackage{relsize}
\usepackage{arydshln}

\newtheorem{e-proposition}[theorem]{Proposition}

\newtheorem{e-definition}[theorem]{Definition\rm}

\renewcommand{\theequation}{\arabic{equation}}
\def\og{\leavevmode\raise.3ex\hbox{$\scriptscriptstyle\langle\!\langle$~}}
\def\fg{\leavevmode\raise.3ex\hbox{~$\!\scriptscriptstyle\,\rangle\!\rangle$}}

\begin{document}

\begin{frontmatter}


\selectlanguage{english}
\title{WENO schemes applied to the quasi-relativistic Vlasov--Maxwell model for laser-plasma interaction}


\selectlanguage{english}
\author[authorlabel1]{Francesco Vecil},
\ead{francesco.vecil@gmail.com}
\author[authorlabel2]{Pep Mulet Mestre},
\ead{jose.mulet@uv.es}
\author[authorlabel3]{Simon Labrunie}
\ead{simon.labrunie@univ-lorraine.fr}

\address[authorlabel1]{Laboratoire de Math\'ematiques, Universit\'e Blaise Pascal (Clermont-Ferrand 2), UMR 6620, CNRS, Campus des C\'ezeaux B.P. 80026, 63171 Aubi\`ere (France)}
\address[authorlabel2]{Universitat de Val\`encia, Departament de Matem\`atica Aplicada, calle del Doctor Moliner 50, Burjassot 46100 (Spain)}
\address[authorlabel3]{  Universit\'e de Lorraine, Institut \'Elie Cartan de Lorraine, UMR 7502, 54506 Vand{\oe}uvre-l\`es-Nancy (France)\\
  CNRS, Institut \'Elie Cartan de Lorraine, UMR 7502, 54506 Vand{\oe}uvre-l\`es-Nancy (France)}

\date{February 7, 2014}

\medskip
\begin{center}
{\small Received *****; accepted after revision +++++\\
Presented by £££££}
\end{center}

\begin{abstract}
In this paper we focus on WENO-based methods for the simulation of the 
1D Quasi-Relativistic Vlasov--Maxwell (QRVM) model used to describe how a laser wave
interacts with and heats a plasma by penetrating into it. 
We propose several non-oscillatory methods based on either Runge--Kutta 
(explicit) or Time-Splitting (implicit) time discretizations.
We then show preliminary numerical experiments.

\vskip 0.5\baselineskip

\selectlanguage{francais}
\noindent{\bf R\'esum\'e}
\vskip 0.5\baselineskip
\noindent
{\bf Sch\'emas WENO appliqu\'es au mod\`ele Vlasov--Maxwell quasi-relativiste pour l'interaction laser-plasma. }
Dans cet article, nous nous int\'eressons aux m\'ethodes de type WENO pour la simulation du mod\`ele
Vlasov--Maxwell quasi-relativiste (QRVM) 1D, utilis\'e pour d\'ecrire la fa\c{c}on dont
une onde laser interagit avec un plasma et le r\'echauffe en le p\'en\'etrant. Nous proposons
plusieurs m\'ethodes non oscillatoires fond\'ees sur
des discr\'etisations en temps soit Runge--Kutta (explicites) soit Time-Splitting (implicites).
Ensuite, nous pr\'esentons des exp\'eriences num\'eriques pr\'eliminaires.

\keyword{Vlasov--Maxwell; WENO; laser-plasma interaction; Runge--Kutta schemes; Strang splitting}
\vskip 0.5\baselineskip
\noindent{\small{\it Mots-cl\'es~:} Vlasov--Maxwell~; WENO~; 
interaction laser-plasma~; sch\'emas de Runge--Kutta~; splitting de Strang}}

\end{abstract}
\end{frontmatter}

\selectlanguage{english}

\section{Introduction}
\noindent The object of our study is the dimensionless 1D quasi-relativistic Vlasov--Maxwell (QRVM) system:
\begin{equation}
\label{QRVM}
\left\{ \begin{array}{l}
\displaystyle \frac{\partial f}{\partial t} + v(p) \frac{\partial f}{\partial x} 
+ F{(t,x)} \frac{\partial f}{\partial p} = 0 \quad \mbox{(collisionless Vlasov)}, 
\qquad F = - \left( E + A B \right)\\[4mm]
\displaystyle  \frac{\partial E}{\partial x} = \eta^{-2} \left( \varrho_{\mathrm{ext}}- \varrho \right) \qquad \mbox{(Poisson equation)}\\[4mm]
\displaystyle \frac{\partial A}{\partial t} = -\mathcal{E}, \qquad
\frac{\partial \mathcal{E}}{\partial t} = \eta^{-2} A\varrho -\frac{\partial B}{\partial x},\qquad
\frac{\partial B}{\partial t} = -\frac{\partial \mathcal{E}}{\partial x},\qquad 
\qquad \mbox{(Maxwell equations)}\\[4mm]
\displaystyle v(p) = \frac{p}{\sqrt{1+p^2}}, \qquad \varrho = \int f \, \mathrm{d}p \qquad \mbox{(relativistic character)}
\end{array} \right.
\end{equation}
solved for $(t,x,p) \in [0,+\infty[ \times [0,1] \times \mathbb{R}$,
endowed with periodic boundary conditions in $x$.
Problem \eqref{QRVM} needs several initial conditions: one for the distribution function $f$; three for the magnetic potential $A$ 
and its derivatives, the magnetic field $B$ and the transverse electric field $\mathcal{E}$, which
are related by
\begin{equation}
\mathcal{E} = -\frac{\partial A}{\partial t} 
\qquad \mbox{and} \qquad
B = \frac{\partial A}{\partial x}.
\label{BE}
\end{equation}
The quantity $\varrho_{\text{ext}}$ represents the immobile ion background which keeps the plasma neutral.
The interest in this Vlasov--Maxwell system is motivated by its importance in plasma physics:
\begingroup
it describes laser-plasma interaction, i.e. the action of a laser wave, 
called \emph{pump}, penetrating into a plasma
and heating it, while interacting with electrostatic waves and accelerating the electrons.
This model, and its variants, have been long known in the plasma physics community~\cite[and references therein]{GBS+90,Huot2003512}.
Its derivation and a discussion about the global existence and uniqueness of classical solutions can be found in \cite{CarLab06}.

\smallbreak

In order to solve~\eqref{QRVM} numerically, one has to choose a time discretization method, a Vlasov solver and a Maxwell solver.
So far, characteristic solvers have been generally used for the Maxwell part, combined with various semi-Lagrangian methods~\cite{GBS+90,Huot2003512,BoCr09} for Vlasov, as well as wavelets~\cite{BLG+08}. Time-splitting methods were often used for the quasi-relativistic model, though they are unstable with a fully relativistic model~\cite{Huot2003512}.

\smallbreak

The goal of this article is to introduce several Weighted Essentially Non-Oscillatory (WENO) schemes for the QRVM model, and to 
perform preliminary tests and comparisons, in order to decide which schemes are more suitable. In Table \ref{overall} we summarize all the combinations we have considered and tested.
\endgroup

\begin{table}[ht]
\centering
\begin{tabular}{|c||l||l|}
\hline
\textbf{time integration}              & RK                            & TS\\\hline
\textbf{Vlasov equation}               & FDWENO                        & \emph{DSLWENO}\\
                                       &                               & \emph{CSLWENO}\\\hline
\textbf{Maxwell equations}             & RK                            & RK\\
                                       & LF                            & LF\\\hline
\end{tabular}
\caption{\textbf{The overall integration strategy.} The schemes in italic are implicit.}
\label{overall}
\end{table}

\noindent \textbf{RK} refers to the Total-Variation-Diminishing \textbf{R}unge--\textbf{K}utta scheme \cite{CCGMS}.\\
\textbf{TS} refers to the \textbf{T}ime-\textbf{S}plitting (Strang) scheme \cite{strang,CarVec05}.\\
\textbf{FDWENO} stands for the \textbf{F}inite-\textbf{D}ifference \textbf{W}eighted \textbf{E}ssentially \textbf{N}on \textbf{O}scillatory 
interpolator for the approximation of partial derivatives \cite{CCGMS}.\\
\textbf{DSLWENO} stands for the non-conservative \textbf{D}irect \textbf{S}emi-\textbf{L}agrangian scheme \cite{CarVec05},
coupled to the Point-Value \textbf{WENO} interpolator \cite{AraBaeBelMul2011,MR2607148,CarVec05}.\\
\textbf{CSLWENO} stands for the \textbf{C}onservative \textbf{S}emi-\textbf{L}agrangian scheme,
based on the Flux-Balance-Method (FBM) \cite{FilSonBer01}, coupled to the FBM\textbf{WENO} described later on.\\
\textbf{LF} stands for \textbf{L}eap-\textbf{F}rog scheme (aka \emph{Yee} scheme).

For the sake of clarity, 
we shall make use of a three-word notation to describe the coupling:
\{\emph{time discretization}\}-\{\emph{Vlasov solver}\}-\{\emph{Maxwell solver}\},
e.g., TS-DSLWENO-LF.

\smallbreak

The outline of this paper is the following: in Section \ref{numericalschemes}
we describe the initial and boundary conditions and the discretization of the system;
in Section \ref{timeintegration} we describe the time-integration strategy;
in Section \ref{numericalresults} we show numerical experiments; 
and we conclude in Section~\ref{weno-is-good-for-you}.

\section{Initial and boundary conditions, and discretization}
\label{numericalschemes}

\subsection{Initialization}
\label{initialize}
\noindent Problem \eqref{QRVM} needs two initializations:
one for the distribution function $f(t,x,p)$, and one for the
electro-magnetic variables $A(t,x)$, $B(t,x)$ and $\mathcal{E}(t,x)$.

\subsubsection{Initialization of the distribution function} 
\noindent We suppose that a proportion $1-\alpha$ of the
electrons are thermalized at a (dimensionless) cold velocity $v_{\mathrm{cold}}$, while
the remaining proportion $\alpha$ are hot with (dimensionless) velocity $v_{\mathrm{hot}}$
\begin{equation*}
G(p) = (1-\alpha) \, G_{\mathrm{cold}}(p) + \alpha \, G_{\mathrm{hot}}(p),
\end{equation*}
where we have split the Maxwellian $G(p)$ into 
a cold part $G_{\mathrm{cold}}(p)$, described by a classical Gaussian,
and a hot part $G_{\mathrm{hot}}(p)$, described by a J\"uttner distribution:
\begin{equation*}
\underbrace{G_{\mathrm{cold}}(p) 
= \frac{\exp\left(-\frac{p^2}{2 \, v_{\mathrm{cold}}^2}\right)}{\sqrt{2 \pi} \, v_{\mathrm{cold}}}}_{\mbox{normalized classical Gaussian}}
\qquad \mbox{and} \qquad
\underbrace{G_{\mathrm{hot}}(p) 
= \frac{\exp\left( -\frac {\sqrt{1+p^2}-1}{v_{\mathrm{hot}}^2} \right)}
{\mathop{\mathlarger{\int}}_{\!\!\!\mathbb{R}} \exp\left( -\frac {\sqrt{1+p^2}-1}{v_{\mathrm{hot}}^2} \right) \, \mathrm{d}p}}_{\mbox{normalized J\"uttner distribution}}.
\end{equation*}

\noindent We shall introduce a fluctuation 
for the initial density
\begin{equation*}
\varrho(0,x) = 1 + \frac{\varepsilon \, v_{\mathrm{cold}} \left( 0.6 \, k_{\mathrm{pla}} \right)}
{ \sqrt{1 + 3 \, v_{\mathrm{cold}}^{2}  \left( 0.6 \, k_{\mathrm{pla}} \right)^{2}} }
\cos\left( 2 \pi \, k_{\mathrm{pla}} \, x \right),
\end{equation*}
for some spatial frequency $k_{\text{pla}}$.
Consequently, a fluctuation is also introduced for the Maxwellian, hence,
all in all, the initial distribution function reads
\begin{equation*}
f(0,x,p) = \varrho(0,x) 
\cdot
G(p - \varepsilon \, v_{\mathrm{cold}} \cos\left( 2 \pi \, k_{\mathrm{pla}} \, x \right)).
\end{equation*}

\subsubsection{Initialization of the electro-magnetic field} 
\noindent 
The initial conditions for $A$, $B$ and $\mathcal{E}$ describe the pump wave
which is going to interact with the plasma wave due to the density fluctuations.

Depending on the coupling we choose between the Vlasov and the Maxwell solvers,
we shall need to set $A$, $B$ and $\mathcal{E}$ at different
initial times and positions, which is why we keep the maximum generality 
by writing them as $(t,x)$-dependent:
\begin{align*}
A\left(t,x\right)           &= A_0 \sin \left(2\pi \, k_{\mathrm{pump}} \, x - \omega_{0} \, t\right), \qquad
B\left(t,x\right)            = 2\pi \, A_0 \, k_{\mathrm{pump}} \cos \left(2\pi \, k_{\mathrm{pump}} \, x - \omega_{0} \, t \right),\\[2mm]
\mathcal{E}\left(t,x\right) &= A_0 \, \omega_{0} \cos \left(2\pi \, k_{\mathrm{pump}} \, x - \omega_{0} \, t\right).
\end{align*}
One checks that the relation between $A$, $B$ and $\mathcal{E}$ at any chosen initial time is given by \eqref{BE}.

\subsubsection{Boundary conditions}
\label{sec_boundcond}
\noindent 
Problem~\eqref{QRVM} is endowed with periodic boundary conditions in the $x$-dimension.
To keep the computational domain bounded and enforce mass conservation, we use 
Neumann boundary conditions in the $p$-dimension.
Actually, if the size of the domain is properly chosen, no
electrons should reach the $p$-border. The boundary conditions are implemented as: 
\begin{align*}
f_{-i,j}    &= f_{i+N_{x},j}      && \mbox{for } i=1,\dots,N_{\mathrm{ghp}}, \quad j=1,\dots,N_{p}-1\\
f_{N_{x}+i,j} &= f_{i           ,j} && \mbox{for } i=1,\dots,N_{\mathrm{ghp}},  \quad j=1,\dots,N_{p}-1\\
f_{i,N_{p}+j} &= f_{i,N_{p}-1}     && \mbox{for } i=0,\dots,N_{x}, \quad j=0,\dots,N_{\mathrm{ghp}}\\
f_{i,-j}    &= f_{i,1}          && \mbox{for } i=0,\dots,N_{x}, \quad j=0,\dots,N_{\mathrm{ghp}}.
\end{align*}

\subsection{Discretization}
\label{amr}
\noindent We mesh the computational domain $\displaystyle \Omega = [0,1] \times \left[-p_{\max}, p_{\max}\right]$
by uniform grids:
\begin{equation*}
\left( x_{i}, p_{j} \right) = \left( i \, \Delta x, j \, \Delta p \right), 
\qquad \!\!\!\! \left( \Delta x, \Delta p \right) = \left( \dfrac{1}{N_{x}}, \dfrac{2 \, p_{\max}}{N_{p}} \right). \nonumber
\end{equation*}
In order to take into account the boundary conditions,
ghost points outside the physical domain are used.

\section{Time integration}
\label{timeintegration}
\noindent In this section, we take care of the time integration 
for the Vlasov equation
\begin{equation}
\frac{\partial f}{\partial t} + v(p) \frac{\partial f}{\partial x} + F(t,x) \frac{\partial f}{\partial p} = 0,
\qquad F = -\left( E + A \,  B \right) 
\label{vlasov}
\end{equation}
and for the set of Maxwell equations
\begin{equation}
\frac{\partial A}{\partial t} = -\mathcal{E}, \qquad
\frac{\partial \mathcal{E}}{\partial t} = \eta^{-2} A\varrho -\frac{\partial B}{\partial x},\qquad
\frac{\partial B}{\partial t} = -\frac{\partial \mathcal{E}}{\partial x}.
\label{Maxwell}
\end{equation}

\noindent As for the Poisson equation
\begin{equation*}
\frac{\partial E}{\partial x} = \eta^{-2} \left( \varrho_{\mathrm{ext}}- \varrho \right),
\end{equation*}
we use the fast, spectrally-accurate solver,
whose details can be found in \cite{MulVec2012}.

\noindent We wish to test two different integration strategies, which are summarized in Table \ref{overall}.

\noindent TS is \emph{implicit} in the sense that it generally uses
implicit schemes for advection, thus weakening the constraints on the time step; 
on the other hand, RK is \emph{explicit}, thus it requires a CFL condition.

\noindent This section is organized as follows:
in Section \ref{RK-schemes} we introduce the Runge--Kutta based schemes;
in Section \ref{TS-schemes} we introduce the Strang-splitting based schemes;
in Section \ref{maxwelleqs} we introduce leap-frog and multi-stage schemes to integrate \eqref{Maxwell};
in Section \ref{checasino!} we summarize all the resulting schemes.

\subsection{RK-FDWENO scheme}
\label{RK-schemes}
\noindent The explicit third-order TVD \emph{Runge--Kutta} strategy consists in
integrating, from time $t^{n}$ to $t^{n+1}$, the Vlasov equation
\begin{equation*}
\frac{\partial f}{\partial t} = - v(p) \frac{\partial f}{\partial x} - F(t,x) \frac{\partial f}{\partial p} =: \mathcal{H}\left[ t, f \right]
\end{equation*}
as
\begin{align}
f^{n,1} &= f^{n} + \Delta t \, \mathcal{H}\left[ t^{n}, f^{n} \right],
\quad
f^{n,2} = \frac{3}{4} f^{n} + \frac{1}{4} f^{n,1} + \frac{1}{4} \Delta t \, \mathcal{H}\left[ t^{n}+\Delta t, f^{n,1} \right], \nonumber\\
f^{n+1} &= \frac{1}{3} f^{n} + \frac{2}{3} f^{n,2} + \frac{2}{3} \Delta t \, \mathcal{H}\left[ t^{n}+\frac{\Delta t}{2}, f^{n,2} \right].
\label{TVD-RK3}
\end{align}

\noindent The partial derivatives are approximated through the fifth-order FDWENO routine for
finite differences, whose details can be found, for instance,
in \cite{CGMS-JCP,CGMS-JCP2} and references therein. 
As this scheme is quite classical, we believe it does not deserve further details here.
The scheme is subject to a CFL constraint for stability:
\begin{equation*}
\Delta t < \dfrac{1}{ \dfrac{\left\| v(p) \right\|_{\infty}}{\Delta x} + \dfrac{\left\| F \right\|_{\infty}}{\Delta p} }.
\end{equation*}
Remark that we have to use the correct upwinding and that,
with proper boundary conditions (see Section~\ref{sec_boundcond}), the scheme enforces mass conservation.

\noindent RK requires the calculation of the Lorentz force 
at three different times 
\begin{equation*}
F\left( t^{n} \right) =: F^{n}, \quad
F\left( t^{n}+\Delta t \right) =: F^{n+1}, \quad
F\left( t^{n}+\frac{1}{2} \Delta t \right) =: F^{n+1/2}.
\end{equation*}
Computing the electrostatic field $E(t)$ at the desired times is easy,
because it is consistent with the distribution function $f(t)$; conversely, obtaining
the magnetic variables $A(t)$ and  $B(t)$ is slightly more complicated, because they follow
their own evolution equations. In case the time integrator for the Maxwell equations
does not provide us with $A$ and $B$ at the desired times,
we can estimate them by interpolations.

\subsection{TS-DSLWENO and TS-CSLWENO schemes}
\label{TS-schemes}
\begin{figure}[ht]
\centering
\subfigure[Prediction.]{\includegraphics[height=4.5cm,keepaspectratio]{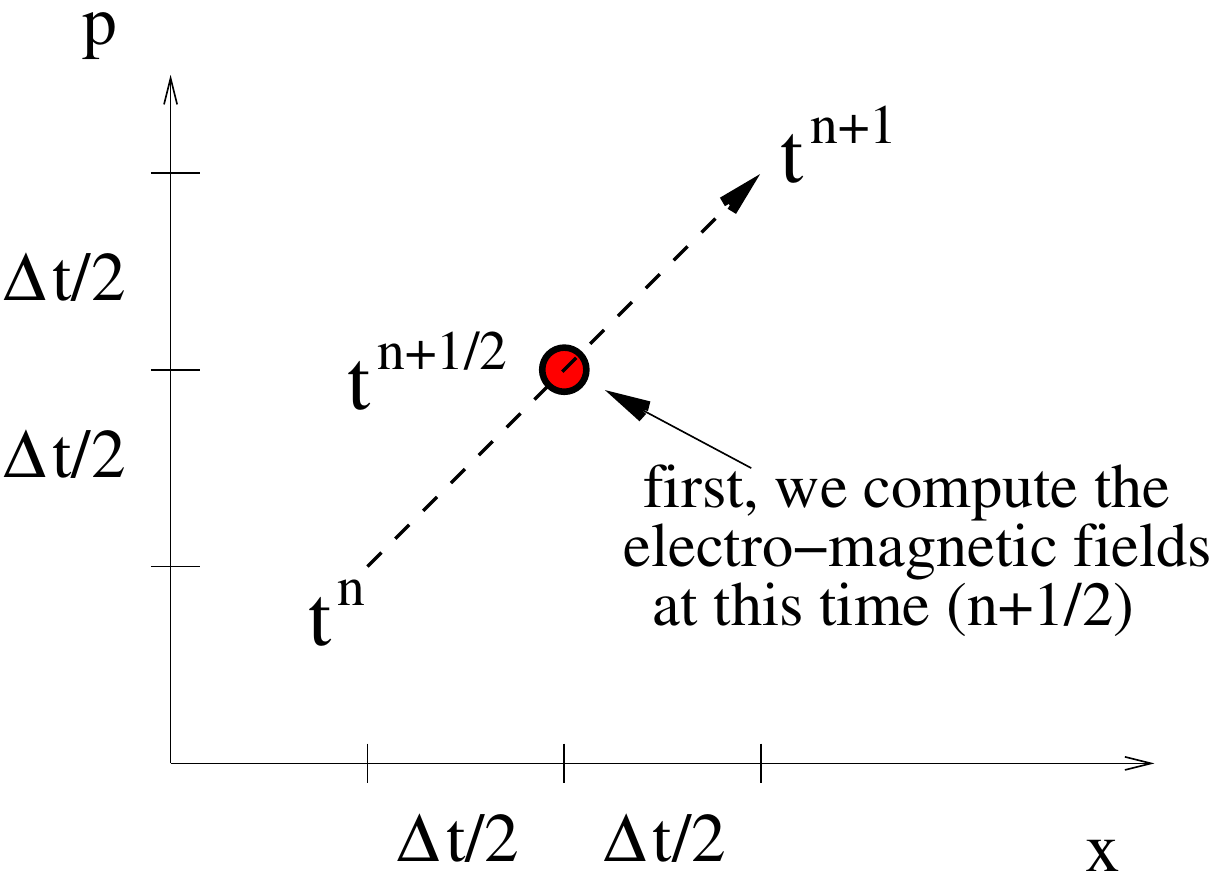}\label{here}}
\hspace{1cm}
\subfigure[Correction.]{\includegraphics[height=4.5cm,keepaspectratio]{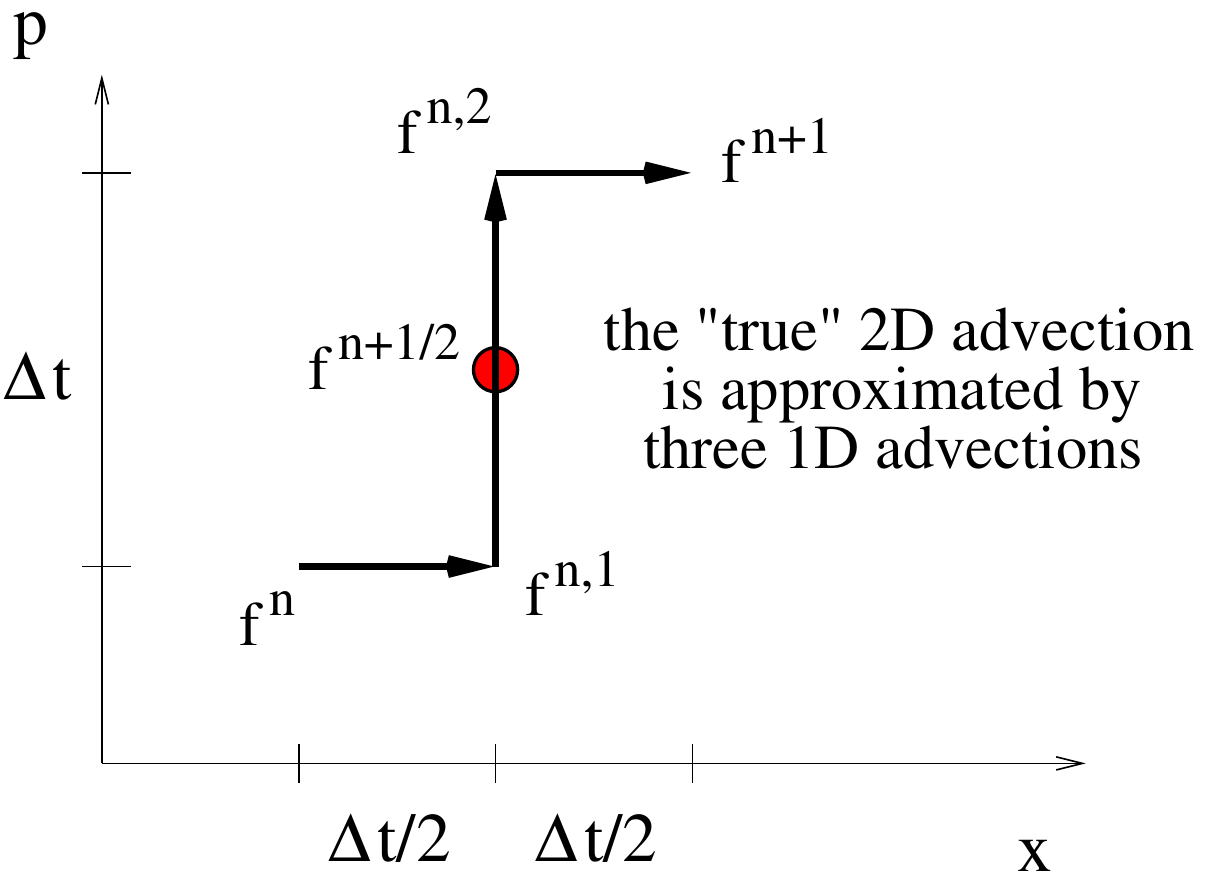}\label{nomodifydensity}}
\caption{\textbf{Prediction/correction strategy.}}
\label{predcorr}
\end{figure}

\noindent The (Strang) \emph{Time-Splitting} strategy \cite{strang,CheKno76} approximates the integration of the Vlasov equation
\begin{equation*}
\frac{\partial f}{\partial t} + v(p) \frac{\partial f}{\partial x} + F{(t,x)} \frac{\partial f}{\partial p} = 0
\end{equation*}
as a combination of partial solutions along the $x$-dimension and the $p$-dimension:
\begin{equation}
\frac{\partial f}{\partial t} + v(p) \frac{\partial f}{\partial x} = 0
\qquad \mbox{and} \qquad
\frac{\partial f}{\partial t} + F{(t,x)} \frac{\partial f}{\partial p} = 0.
\label{both}
\end{equation}
We advect $f^{n} \longmapsto f^{n+1}$ by means of the advection 
field evaluated at time $t^{n+1/2}$, a strategy called \emph{prediction/correction} \cite{JakMarOwr2000,MarOwr2001},
summarized on Figure \ref{predcorr},
which gives a scheme of order~2 in time as soon as $F(t^{n+1/2})$ is approximated at order~1. 

\noindent In principle, the one-dimensional PDEs \eqref{both} can be solved by means of any time integrator; 
here we propose a direct semi-Lagrangian (DSL) strategy (non-conservative),
fully described in \cite{MulVec2012},
and a conservative semi-Lagrangian (CSL) strategy,
described in Section \ref{cslweno};
\emph{semi-Lagrangian} means that the method is characteristics-based.

\subsubsection{CSL integration for 1D advection problems}
\label{cslweno}
\noindent The model equation which we solve is
\begin{equation*}
\frac{\partial u}{\partial t} + \frac{\partial}{\partial x} \left[ a(t,x) \, u \right] = 0, 
\qquad u(t^{\star},x) = u^{\star}(x), \qquad (t,x) \in [0,+\infty[ \times I
\end{equation*}
(being $a : [0,+\infty[ \times I \to \mathbb{R}$ and $I \subseteq \mathbb{R}$ an interval)
by means of a semi-Lagrangian conservative method; this strategy
is taken from \cite{FilSonBer01}.
To this end, we evolve approximated cell averages
\begin{equation*}
u^{n+1}_{i} \approx {\frac1{\Delta x}\int_{x_{i-1/2}}^{x_{i+1/2}} } u\left( t^{n+1}, \xi \right) \, \mathrm{d}\xi
\end{equation*}
and use a semi-Lagrangian strategy by following the characteristics backward,
along which $J \, u$ is conserved,
\begin{equation}
\int_{x_{i-1/2}}^{x_{i+1/2}} u\left( t^{n+1}, \xi \right) \, \mathrm{d}\xi
= \int_{x_{i-1/2}}^{x_{i+1/2}} u\left( t^{n}, \mathcal{X}\left( t^{n}; t^{n+1}, \xi \right) \right) \, 
J\left( t^{n}; t^{n+1}, \xi \right) \, \mathrm{d}\xi,
\label{fbm}
\end{equation}
with $\mathcal{X}(s;t,x)$ the characteristic and $J(s; t, x)$ its Jacobian:
\begin{equation*}
\frac{\mathrm{d} \mathcal{X}(s; t, x)}{\mathrm{d}s} = a\left(s, \mathcal{X}(s; t, x)\right), \qquad \mathcal{X}(t;t,x) = x,
\qquad J(s; t, x) := \det \dfrac{\partial \mathcal{X}(s; t, x)}{\partial x}.
\end{equation*}
If we change variables $\displaystyle \eta = \mathcal{X}\left( t^{n}; t^{n+1}, \xi \right)$
into \eqref{fbm}, we get:
\begin{equation}
{\frac1{\Delta x}\int_{x_{i-1/2}}^{x_{i+1/2}} } u\left( t^{n+1}, \xi \right) \, \mathrm{d}\xi
= \frac{1}{\Delta x} \int_{x^{\mathrm{back}}_{i-1/2}}^{x^{\mathrm{back}}_{i+1/2}}
u\left( t^{n}, \eta \right) \, \mathrm{d}\eta
= \frac{U^{n}\left( x^{\mathrm{back}}_{i+1/2} \right) - U^{n}\left( x^{\mathrm{back}}_{i-1/2} \right)}{\Delta x},
\label{fbmscheme0}
\end{equation}
where we have set $\displaystyle x^{\mathrm{back}} := \mathcal{X}\left( t^{n}; t^{n+1}, x \right)$
and  $U^{n}$ is a primitive of $u(t^{n}, \cdot)$. This gives the following scheme:
\begin{equation}
u_{i}^{n+1}
= \frac{\widetilde U^{n}\left( x^{\mathrm{back}}_{i+1/2} \right) -
  \widetilde U^{n}\left( x^{\mathrm{back}}_{i-1/2} \right)}{\Delta x},
\label{fbmscheme}
\end{equation}
where $\widetilde U^{n}$ is an approximation of $U^{n}$ based on values of $\left( u_{j}^{n} \right)_j$.
The scheme is conservative if $u$ is compactly supported or under periodic boundary conditions.
In our application, the computations are simplified by $a$ being a real constant.
\footnote{Recall that the advection field in the $x$-dimension is independent of~$x$, and similarly in the $p$-dimension; furthermore $F(t,x)$ is approximated by $F(t^{n+1/2},x)$ on the time interval $[t^n,t^{n+1}]$.}
Therefore, we have explicit characteristics $\mathcal{X}(s; t, x) = x + a(s-t)$, so
\begin{equation*}
u^{n+1}_{i} 
= \frac{\widetilde U^{n}\left( x_{i+1/2} - a \, \Delta t \right) -
  \widetilde U^{n}\left( x_{i-1/2} - a \, \Delta t \right)}{\Delta x}.
\end{equation*}

\subsubsection{The WENO reconstruction for CSL (called FBMWENO)}
\label{cslweno2}
\noindent In order to set up the scheme \eqref{fbmscheme} we need an interpolator for the primitive $U$
(dropping the time-dependency notation from now on).
In the WENO fashion, we shall perform a convex combination of several Lagrange polynomials
interpolating $U$ at different substencils.
We can adjust two parameters in order to obtain all the possible combinations: 
the degree $r_{\mathrm{tot}}$ of the Lagrange polynomial interpolating $U(x)$ in the whole stencil $\mathcal{S}$
(which thus contains $r_{\mathrm{tot}}+1$ points), 
and the degree $r_{\mathrm{sub}}$ of the Lagrange polynomials
in the substencils (each substencil contains $r_{\mathrm{sub}}+1$ points). 
Let us also introduce the number of substencils $N_{\mathrm{sub}} := r_{\mathrm{tot}}-r_{\mathrm{sub}}+1$.

\noindent 
Let us denote $P_{\nu}^{r}(x)$ the Lagrange polynomial interpolating the point values of the primitive $U$
at points $\left\{ x_{\nu-r}, \dots, x_{\nu} \right\}$.
If $\displaystyle \mathcal{S} = \left\{ x_{\mathrm{left}} , \dots, x_{\mathrm{left}+r_{\mathrm{tot}}} =: x_{\mathrm{right}} \right\}$ is
the big stencil used to approximate $U(x)$, then
\begin{equation*}
U(x) \approx \widetilde U(x):=\sum_{\ell = 0}^{N_{\mathrm{sub}}-1} \omega_{\ell}(x) \, P^{r_{\mathrm{sub}}}_{\mathrm{right}-\ell}(x).
\end{equation*}
In order to define the weights 
\begin{equation*}
\omega_{\ell}(x) := \frac{\tilde\omega_{\ell}(x)}{\sum_{\ell'=0}^{N_{\mathrm{sub}}-1} \tilde\omega_{\ell'}(x)}, \qquad 
\tilde\omega_{\ell}(x) := \frac{C_{\ell}(x)}{\left( 10^{-6} + \sigma_{\ell} \right)^{2}}, \qquad \ell=0, \dots, N_{\mathrm{sub}}-1
\end{equation*}
we need two ingredients: the polynomials $\left\{ C_{\ell}(x) \right\}_{\ell=0}^{N_{\mathrm{sub}}-1}$
defined by the relation 
\begin{equation*}
P_{\mathrm{right}}^{r_{\mathrm{tot}}}(x) = \sum_{\ell=0}^{N_{\mathrm{sub}}-1} C_{\ell}(x) \, P_{\mathrm{right} - \ell}^{r_{\mathrm{sub}}}(x),
\end{equation*}
and the \emph{smoothness indicators} $\left\{ \sigma_{\ell} \right\}_{\ell=0}^{N_{\mathrm{sub}}-1}$,
which we wish to define in such a way that the $u^{n+1}_{i}$ given by~\eqref{fbmscheme}
is not polluted by spurious oscillations.
To this end, we are not interested in the smoothness of $U$, rather in the smoothness of $u$.

\noindent Now, the derivative of $P^r_{\nu}(x)$ is a lower-order approximation to $u(x)$:
\begin{equation*}
\mathcal{P}^r_{\nu}(x) := \frac{\mathrm{d} P^r_{\nu}}{\mathrm{d} x} (x) \approx u(x),
\end{equation*}
in the sense that if $P^r_{\nu}(x)$ approximates $U(x)$ at order $r+1$,
$\mathcal{P}^r_{\nu}(x)$ approximates $u(x)$ at order $r$.
We now fix the interval $I := \left[ x_{i-1/2}, x_{i+1/2} \right]$
that contains the evaluation point and define the smoothness measurement as
in the Jiang--Shu fashion \cite{JS}: for $\ell = 0, \dots, N_{\mathrm{sub}}-1$
\begin{equation*}
\sigma_{\ell} := 
\sum_{k=1}^{r_{\mathrm{sub}}-1} \Delta x^{2k-1} \int_{I} \left[ \left( \mathcal{P}_{\mathrm{right} - \ell}^{r_{\mathrm{sub}}} \right)^{(k)} (\xi) \right]^{2} \, \mathrm{d}\xi.
\end{equation*}
\noindent The polynomials $C_\ell(x)$ and constants $\sigma_\ell$ for $\left( r_{\mathrm{tot}}, r_{\mathrm{sub}}
\right) = (5,3)$ are given in Appendix \ref{constants_fbmweno}.

\subsection{Integration of the Maxwell equations}
\label{maxwelleqs}
\noindent We test two strategies: a \emph{leap-frog}-type Yee scheme
and a Runge--Kutta scheme. The Yee scheme will be coupled to both
schemes for the Vlasov equation and the Runge--Kutta scheme will only
be coupled to the Runge--Kutta scheme for the Vlasov equation.

In any case, once we have updated the ponderomotive force
$\displaystyle \mathcal{F} = A \, B$ 
up to time $t^{n}$, we impose it has numerically zero average:
\begin{equation*}
\bar{\mathcal{F}}^{n} := \frac{1}{N_{x}} \sum_{i=0}^{N_{x}-1} \mathcal{F}^{n}_{i},
\qquad \mbox{then} \qquad \mathcal{F}^{n}_{i} \longmapsto \mathcal{F}^{n}_{i} - \bar{\mathcal{F}}^{n}.
\end{equation*}

\begin{figure}
\centering
\subfigure[leap $A^{n}_{i} \mapsto A^{n+1}_{i}$]{\includegraphics[scale=.4]{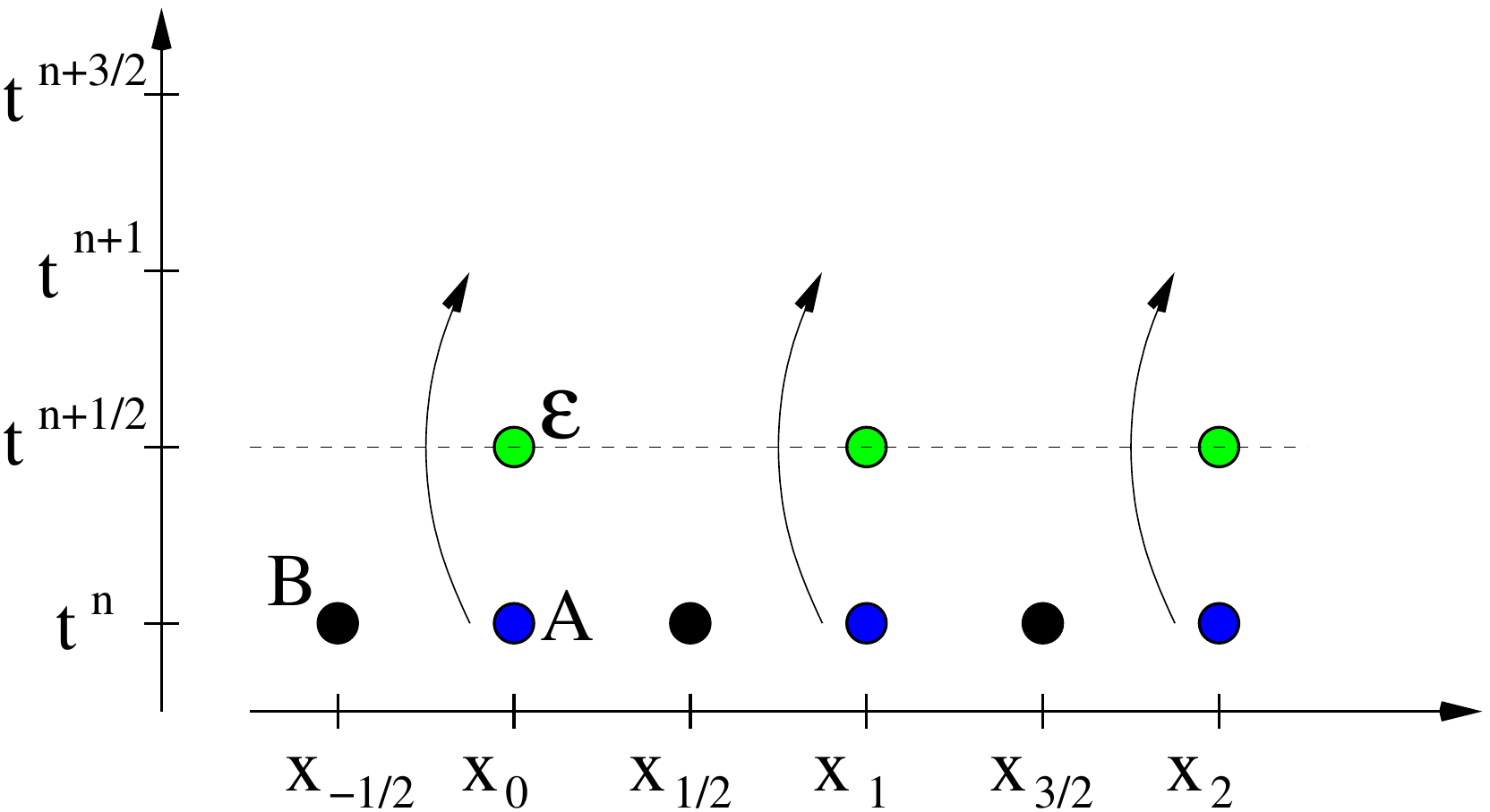}\label{yee_1}}
\subfigure[leap $B^{n}_{i+1/2} \mapsto B^{n+1}_{i+1/2}$]{\includegraphics[scale=.4]{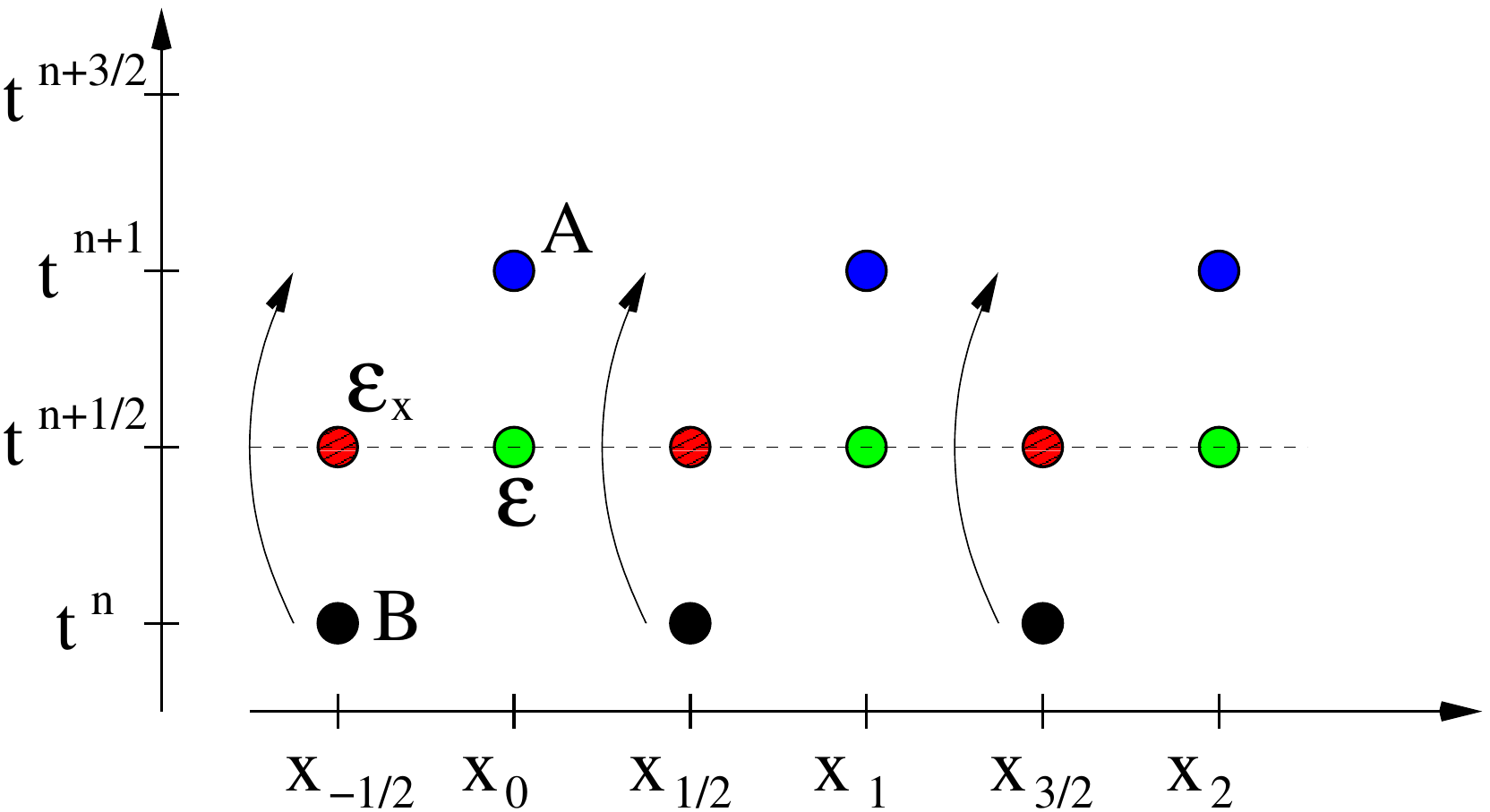}\label{yee_2}}
\subfigure[leap $\mathcal{E}^{n+1/2}_{i} \mapsto \mathcal{E}^{n+3/2}_{i}$]{\includegraphics[scale=.4]{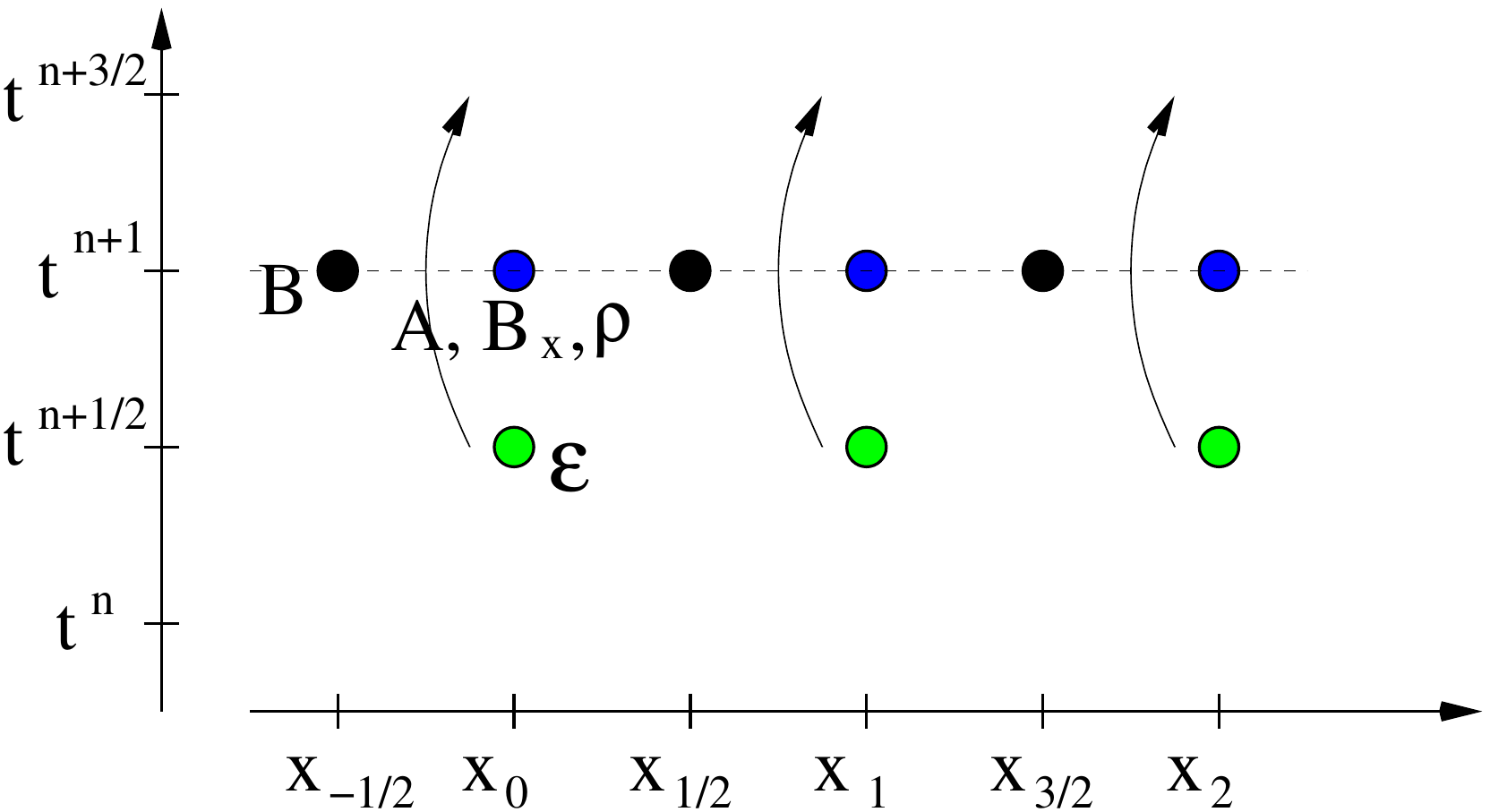}\label{yee_3}}
\subfigure[at time $t^{n+1}$]{\includegraphics[scale=.4]{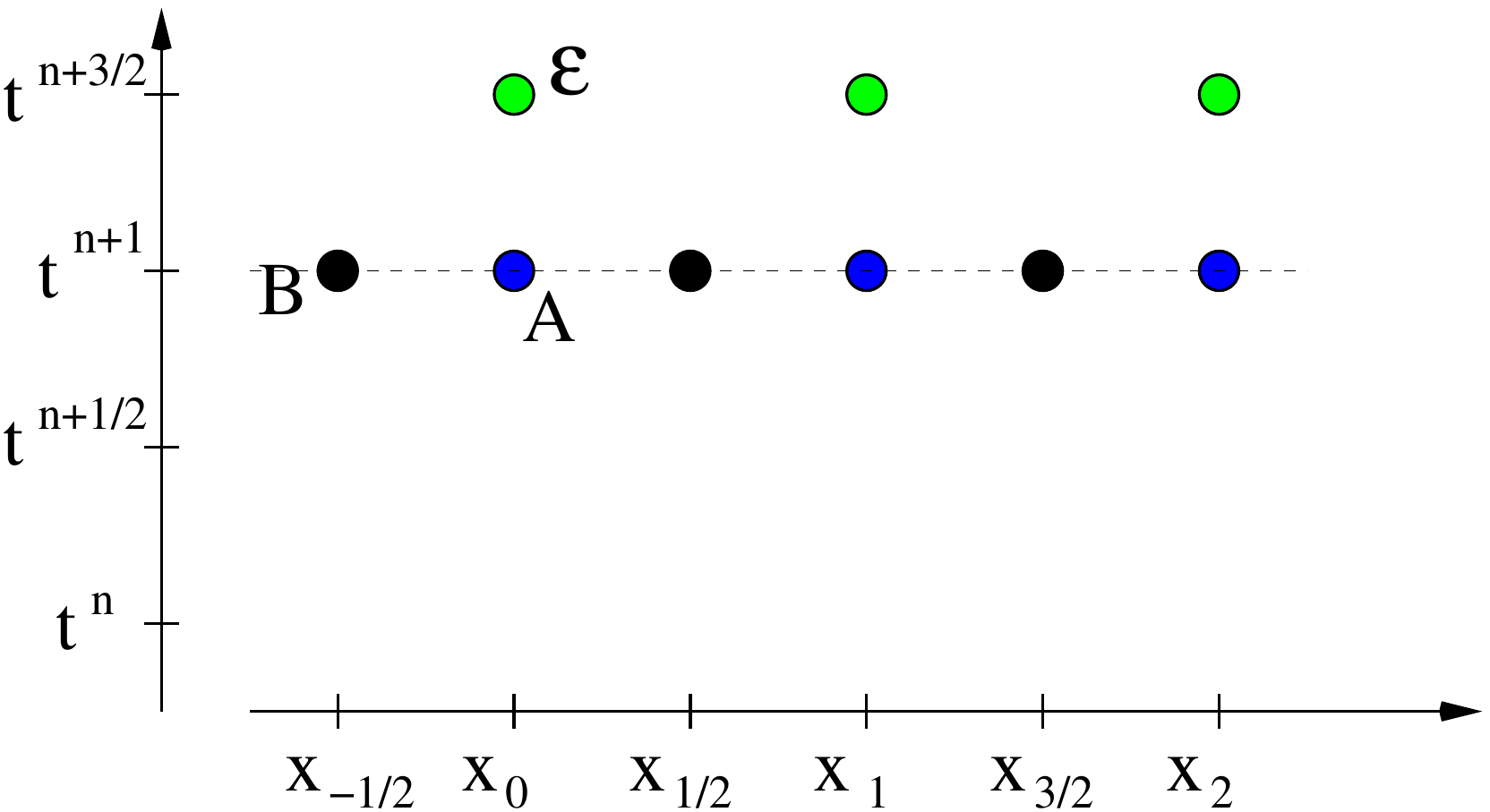}\label{yee_4}}
\caption{\textbf{\emph{Leap-frog} strategy.} The scheme is second-order in both time and space,
  because all the $t$- and $x$-derivatives are approximated by centered differences.
  Inside the figure $\displaystyle \mathcal{E}_{x} := \frac{\partial \mathcal{E}}{\partial x}$ 
  and $\displaystyle  B_{x} := \frac{\partial B}{\partial x}$.}
\label{leapfrog}
\end{figure}

The LF scheme that we use for the Maxwell equations is second-order accurate in both space and time,
and is known as the \emph{Yee} scheme. It is of the \emph{leap-frog} type
with half-shifted variables: see Figure \ref{leapfrog} for a sketch.
Knowing $\varrho^{n+1}$, we advance in time 
$\left( A^{n}, B^{n}, \mathcal{E}^{n+1/2} \right) \longmapsto \left( A^{n+1}, B^{n+1}, \mathcal{E}^{n+3/2} \right)$ 
by centered finite differences:

\begin{itemize}
\item The evolution of the vector potential $A$ (Figure \ref{yee_1})
\begin{equation*}
\frac{\partial A}{\partial t} = -\mathcal{E}
\qquad \mbox{gives} \qquad 
A^{n+1}_{i}=A^{n}_{i} -\mathcal{E}^{n+1/2}_{i} \Delta t.
\end{equation*}
\item The evolution of the magnetic field $B$ (Figure \ref{yee_2}), 
\begin{equation*}
\frac{\partial B}{\partial t} = -\frac{\partial \mathcal{E}}{\partial x}
\qquad
\mbox{gives}
\qquad
B^{n+1}_{i+1/2}=B^{n}_{i+1/2}-\frac{\Delta
t}{\Delta x}\left( \mathcal{E}^{n+1/2}_{i+1}-\mathcal{E}^{n+1/2}_{i} \right).
\end{equation*}
\item The evolution of the transverse electric field $\mathcal{E}$ (Figure \ref{yee_3}),
\begin{equation*}
\frac{\partial \mathcal{E}}{\partial t}
= \eta^{-2} A \, \varrho -\frac{\partial B}{\partial x} 
\quad
\mbox{gives}
\quad
\mathcal{E}^{n+3/2}_{i} =
\mathcal{E}^{n+1/2}_{i}+\eta^{-2} A^{n+1}_{i} \varrho^{n+1}_{i} \Delta t -
\frac{\Delta t}{\Delta x} \left( B^{n+1}_{i+1/2}-B^{n+1}_{i-1/2} \right).
\end{equation*}
\end{itemize}

\subsection{Summary of the schemes}
\label{checasino!}
\noindent In order to construct the schemes resulting from the different choices for the time integrators of
the Vlasov and the Maxwell equations (see Table \ref{overall}), we have to be particularly careful 
in order to fit each block properly within the coupling.

\subsubsection{TS-DSLWENO-LF and TS-CSLWENO-LF schemes}
\begin{figure}[ht]
\centering
\subfigure[Strang 1/3]{\includegraphics[scale=.34]{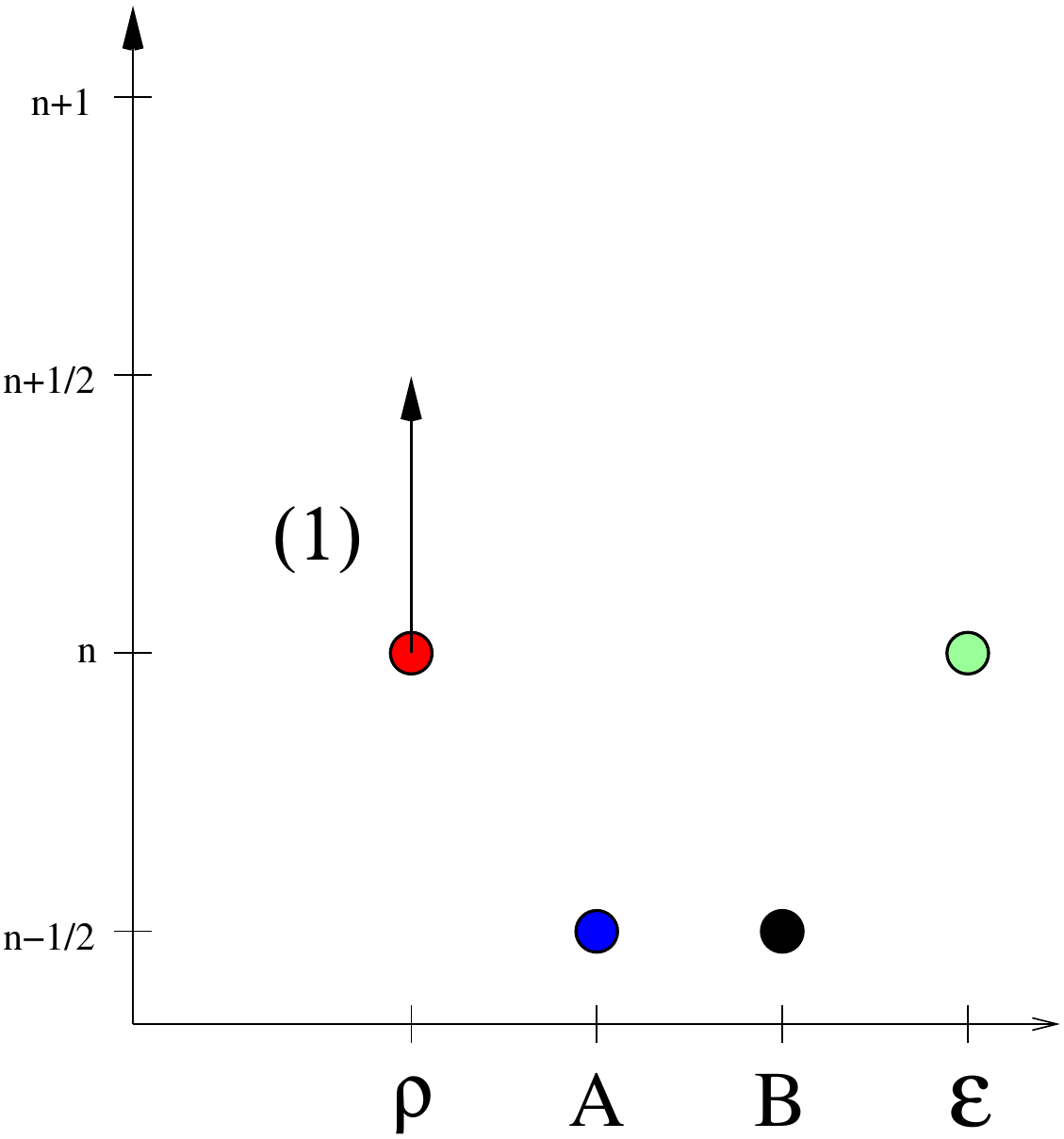}\label{init_TS_LF_1.pdf}}
\subfigure[leap-frog]{\includegraphics[scale=.34]{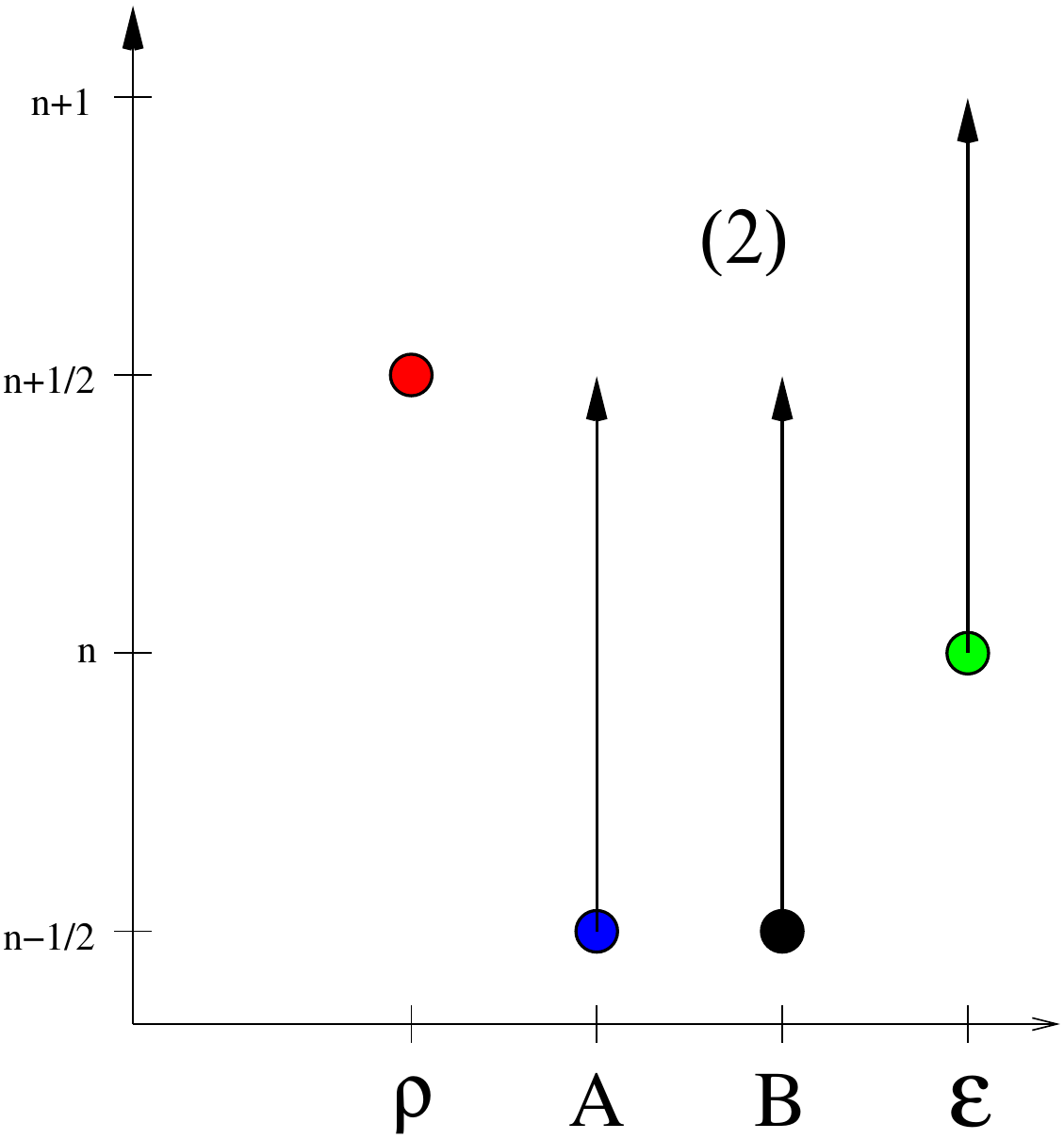}\label{init_TS_LF_2.pdf}}
\subfigure[Strang 2/3, 3/3]{\includegraphics[scale=.34]{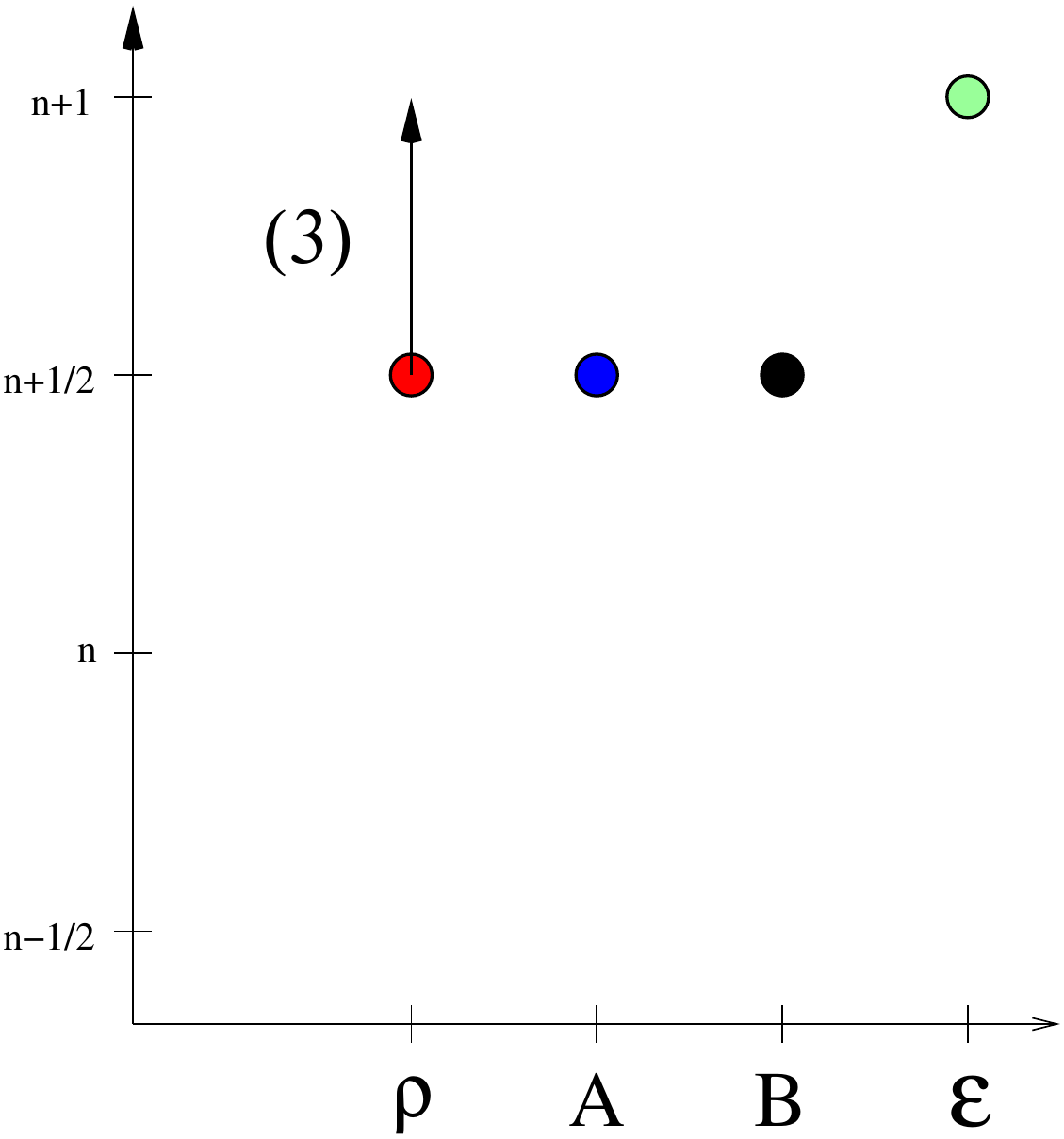}\label{init_TS_LF_3.pdf}}
\subfigure[at time $t^{n+1}$]{\includegraphics[scale=.34]{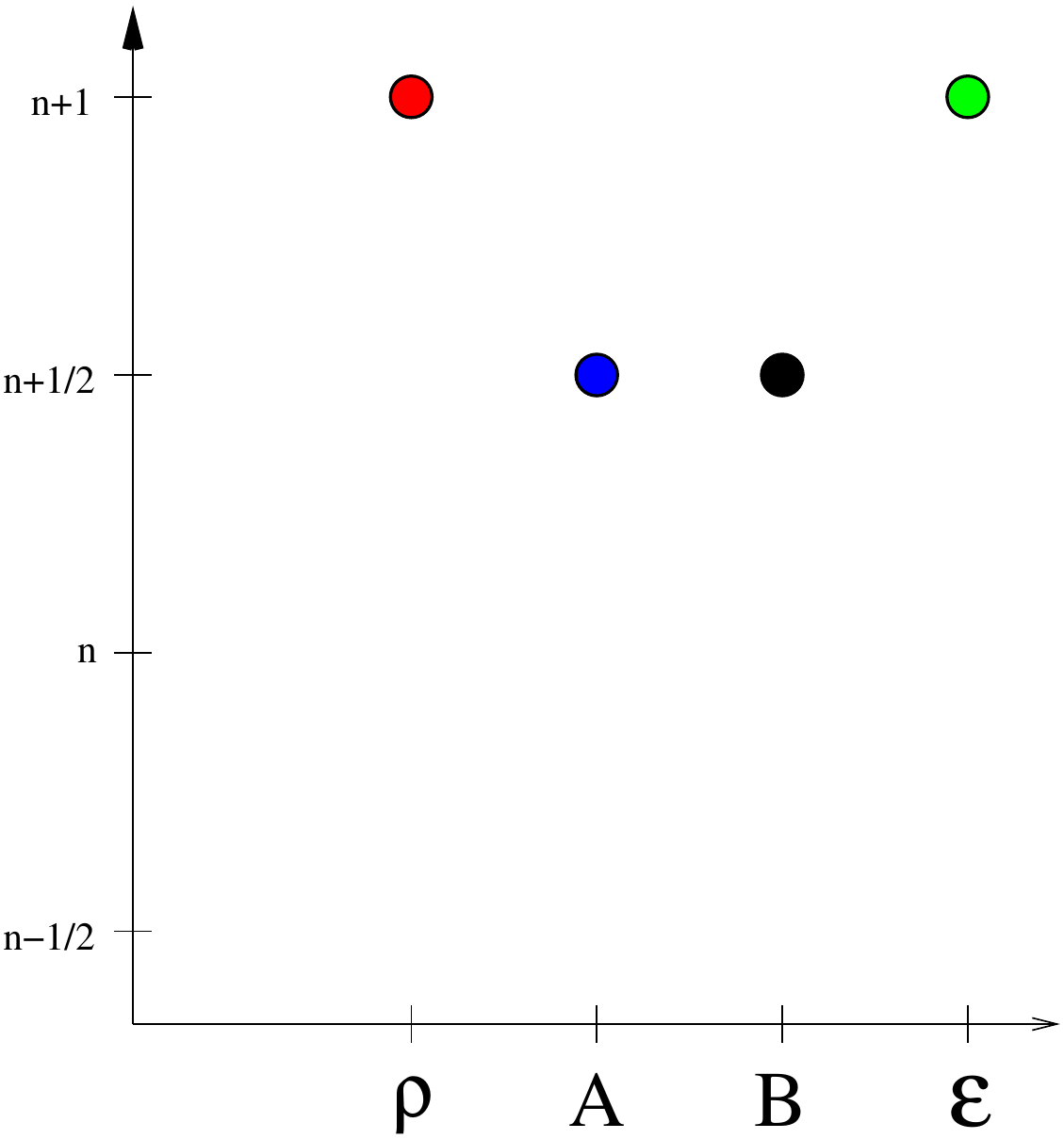}\label{init_TS_LF_4.pdf}}
\caption{\textbf{TS-DSLWENO-LF and TS-CSLWENO-LF schemes.} The schemes differ in how the Strang stages are performed.}
\label{TS-DSLWENO-LF}
\end{figure}

\noindent The scheme  to advance
\begin{equation*}
\left( f^{n}, A^{n-1/2}, B^{n-1/2}, \mathcal{E}^{n} \right) \longmapsto \left( f^{n+1}, A^{n+1/2}, B^{n+1/2}, \mathcal{E}^{n+1} \right)
\end{equation*}
is sketched on Figure \ref{TS-DSLWENO-LF}. 
Notice that the time indices of $A,\ B,\ \mathcal{E}$ have been shifted by one half w.r.t. Section~\ref{maxwelleqs}, so as to have the force at hand at time~$t^{n+1/2}$, as explained in Section~\ref{TS-schemes}. Thus, $\varrho$ and~$E$ must be available at~$t^{n+1/2}$. This is done by computing them after the first half-advection in~$x$~\cite{CheKno76}, see Figure~\ref{nomodifydensity}.
The difference between the two schemes is
how the steps in Figure \ref{init_TS_LF_1.pdf} and Figure \ref{init_TS_LF_3.pdf} are performed,
with a non-conservative method for DSLWENO and with a conservative one for CSLWENO.

\subsubsection{RK-FDWENO-RK scheme}
This scheme is obtained by applying the third-order TVD
Runge--Kutta ODE solver \eqref{TVD-RK3} to a  discretization in $x$ and $p$ of the
Vlasov--Maxwell equations
\begin{equation*}
\frac{\partial}{\partial t} \left( \begin{array}{c} f \\\hdashline A \\ B \\ \mathcal{E} \end{array} \right)
=
\left( \begin{array}{c} - v(p) \frac{\partial f}{\partial x} - (E+A \, B)(t, x) \frac{\partial f}{\partial p} \\\hdashline -\mathcal{E} \\ -\frac{\partial \mathcal{E}}{\partial x} \\ \eta^{-2} A \, \varrho -\frac{\partial B}{\partial x} \end{array} \right)
=: \mathcal{H} \left[ t, \left( \begin{array}{c} f \\\hdashline A \\ B \\ \mathcal{E} \end{array} \right) \right],
\end{equation*}
where, as mentioned in Section \ref{RK-schemes}, the $x$ and $p$ derivatives in the Vlasov
equation are discretized by WENO finite differences, the $x$ derivatives in
the Maxwell equations are discretized by linear finite differences;
$\varrho$ is discretized by the midpoint quadrature rule,
and $E$ is computed by the Poisson solver.

\subsubsection{RK-FDWENO-LF scheme}
\begin{figure}
\centering
\subfigure[RK 1/3]{\includegraphics[scale=.33]{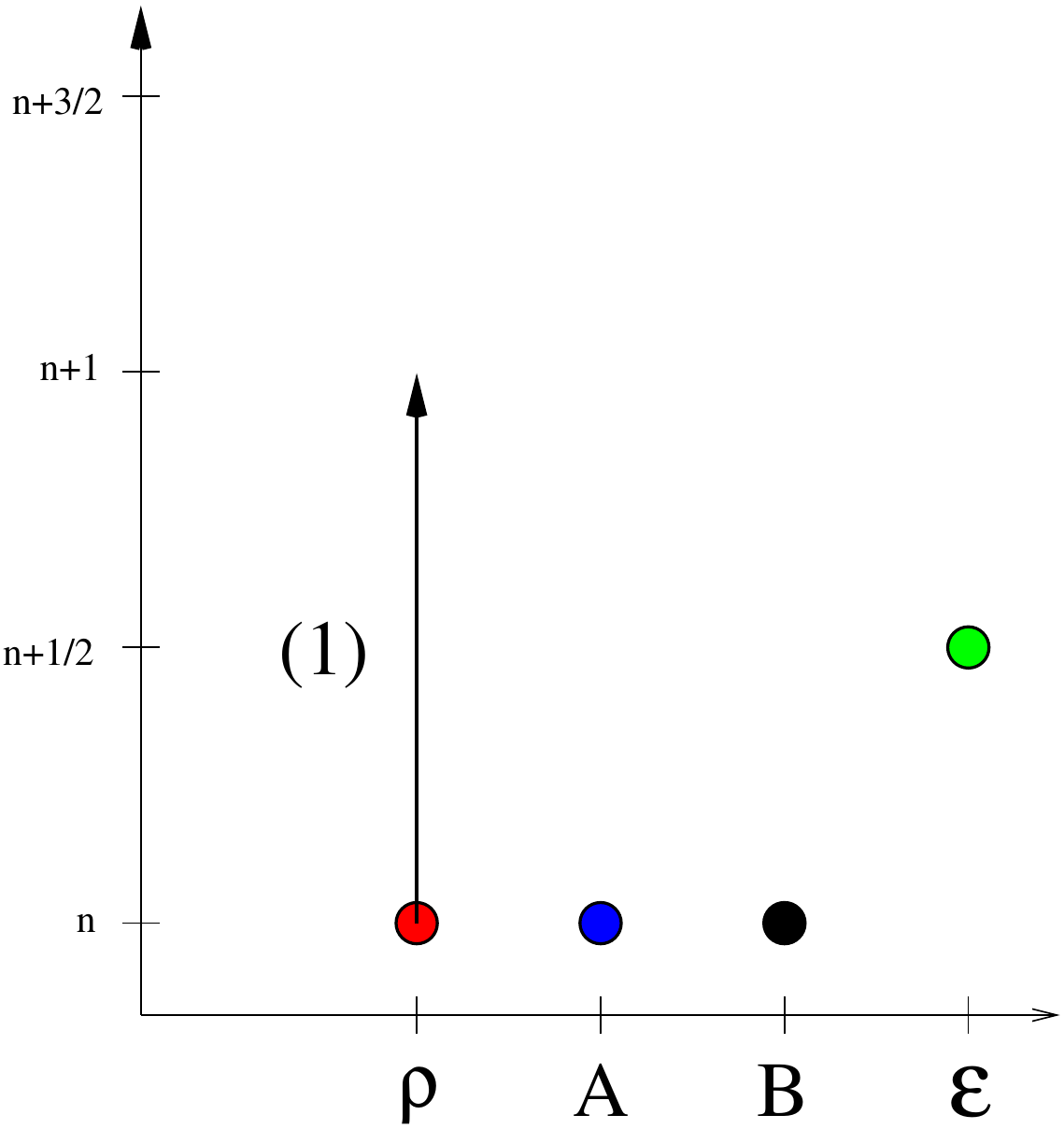}\label{init_RK_LF_1.pdf}}
\subfigure[leap-frog]{\includegraphics[scale=.33]{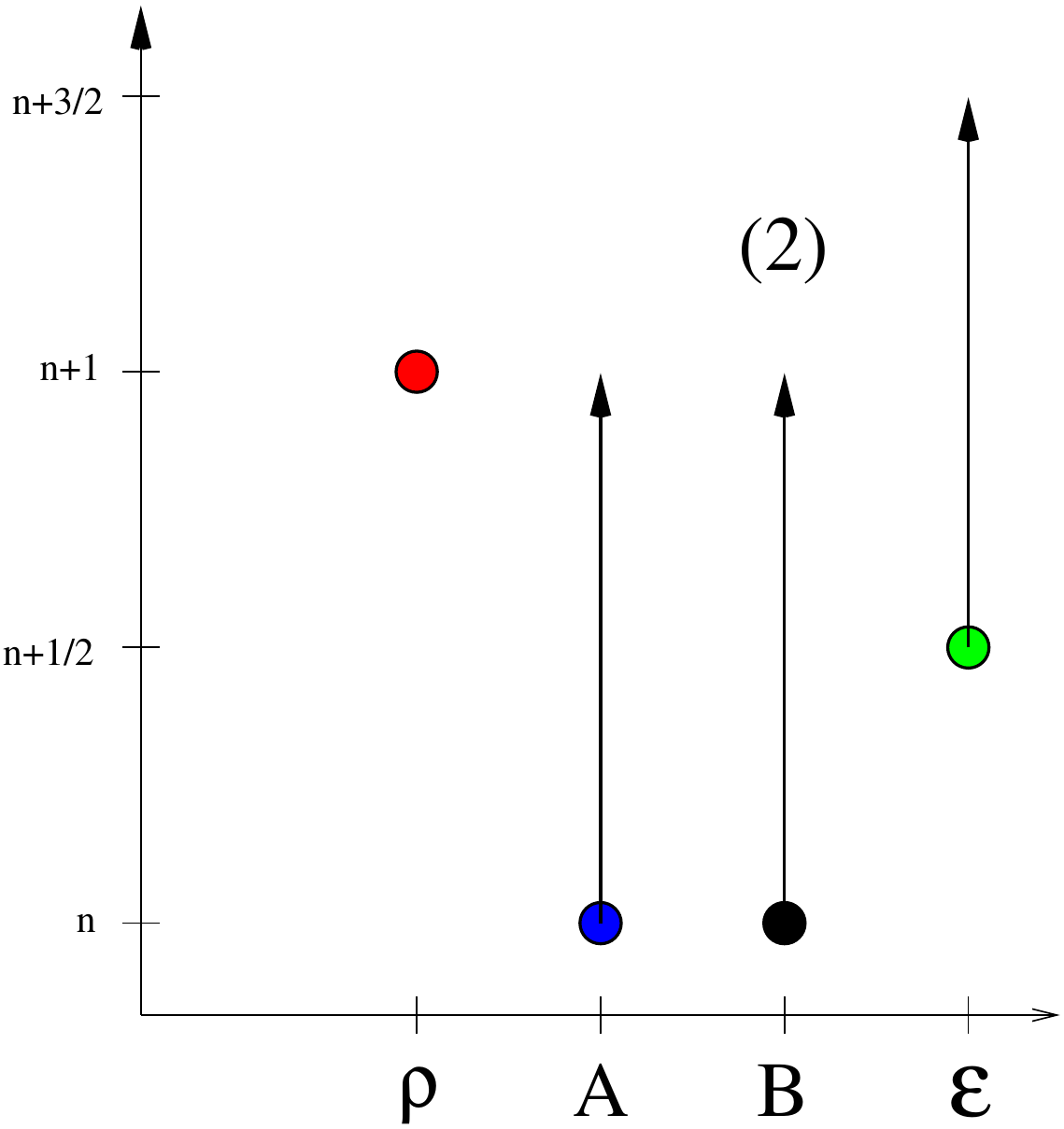}\label{init_RK_LF_2.pdf}}
\subfigure[RK 2/3]{\includegraphics[scale=.33]{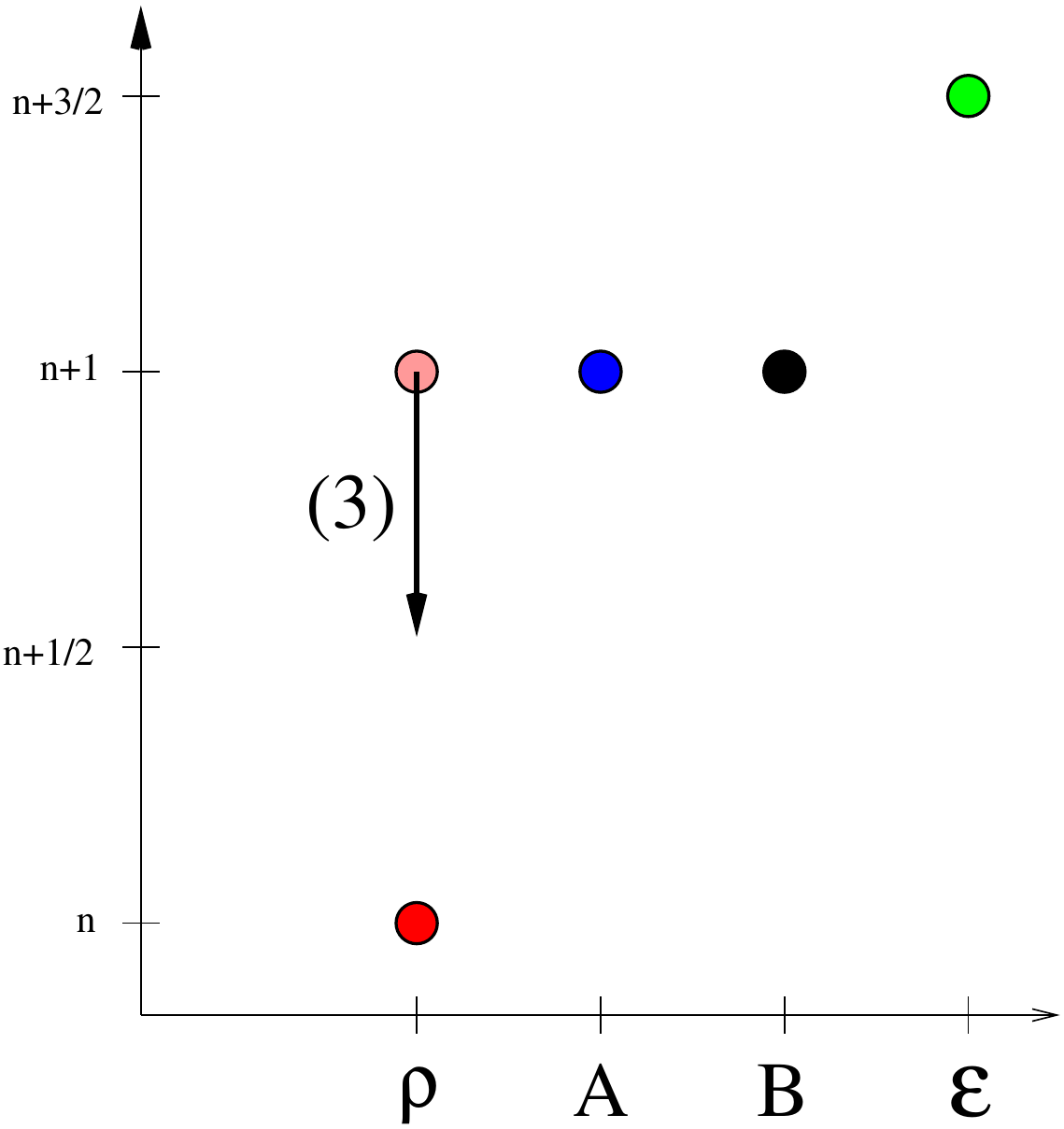}\label{init_RK_LF_3.pdf}}
\subfigure[RK 3/3]{\includegraphics[scale=.33]{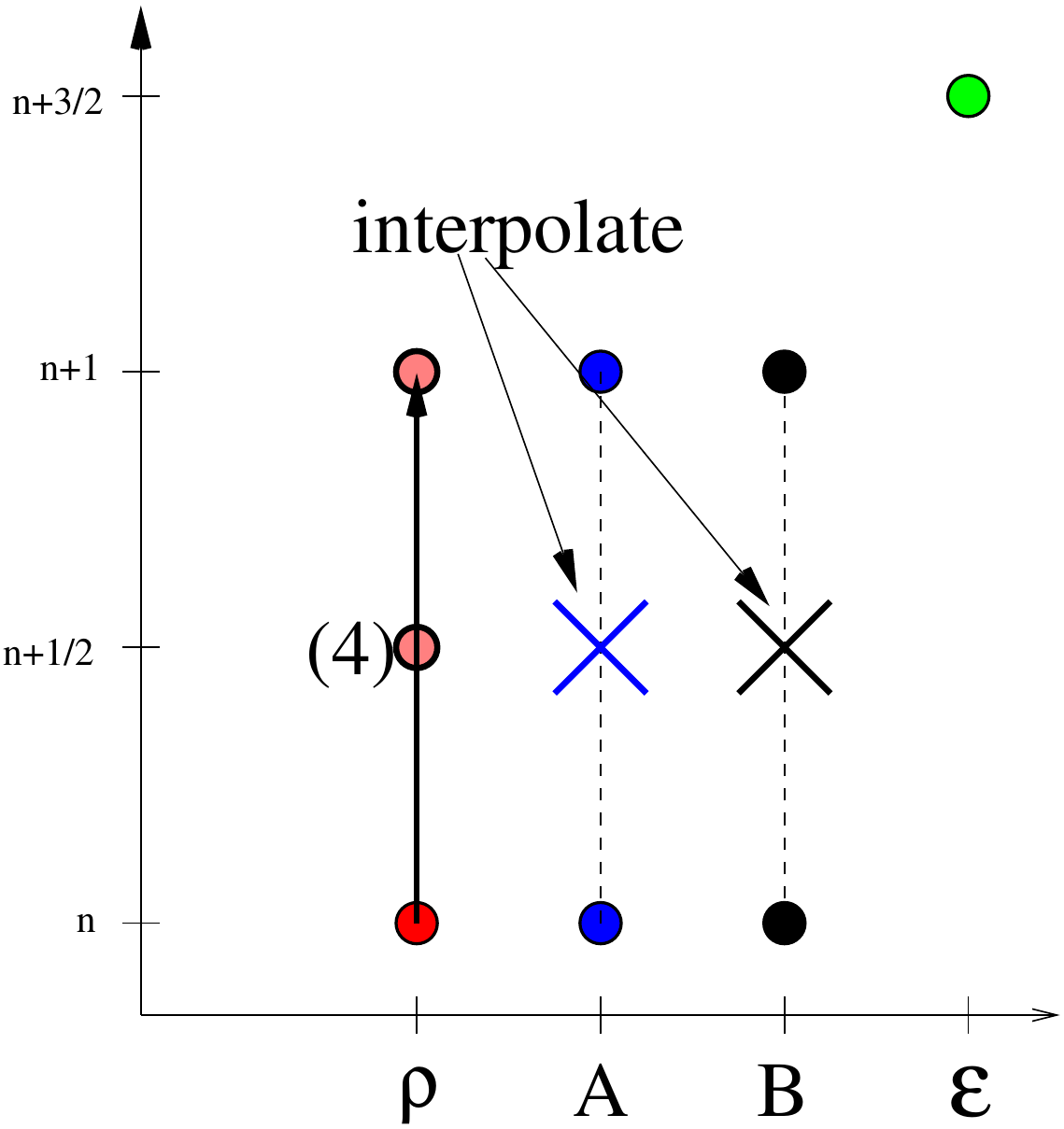}\label{init_RK_LF_4.pdf}}
\caption{\textbf{RK-FDWENO-LF scheme.} This scheme is second-order in time (because of
the interpolation and the first-order approximation of $\varrho$ used to evolve $\mathcal{E}$) and second-order in space.}
\label{RK-FDWENO-LF}
\end{figure}

\noindent The resulting scheme is depicted in Figure
\ref{RK-FDWENO-LF}.
Remark that the Yee scheme forces the time step $\Delta t$ to be kept fixed, 
despite the adaptive character of the Runge--Kutta scheme.

\section{Results for the quasi-relativistic Vlasov--Maxwell system}
\label{numericalresults}
\noindent No WENO-based scheme has yet been extensively tested on the QRVM problem. Therefore,
our first task is to decide which among the overall integration strategies introduced in Table~\ref{overall}
are suitable.

\subsection{Empirical stability results}
\noindent All the schemes proposed in this article seem stable
from empirical observation, but RK-FDWENO-LF requires
extremely small time steps in order not to blow up.
A summary is given in Table \ref{estabilidad}.

\begin{table}[ht]
\centering
\begin{small}
\begin{tabular}{|l||c|c|}
\hline
Vlasov $\downarrow$ \hspace{3cm} Maxwell $\rightarrow$  & \textbf{LF}                                   & \textbf{RK}\\\hline\hline
\textbf{RK}-\textbf{FDWENO}-                            & \scalebox{3}{\frownie}                        & \scalebox{3}{\smiley}\\\hline
\textbf{TS}-\textbf{DSLWENO}-                           & \scalebox{3}{\frownie}                        & not couplable\\\hline
\textbf{TS}-\textbf{CSLWENO}-                           & \scalebox{3}{\smiley}                         & not couplable\\\hline
\end{tabular}
\end{small}
\caption{\textbf{Quality of the results.}}
\label{estabilidad}
\end{table}

\noindent 
The evolution equations for $B$ and $\mathcal{E}$ can be rewritten as
\begin{equation*}
\frac{\partial \left( B \pm \mathcal{E} \right)}{\partial t} 
\pm \frac{\partial \left( B \pm \mathcal{E} \right)}{\partial x}
= \pm \eta^{-2} \, A \, \varrho,
\end{equation*}
therefore the condition 
\begin{equation}
\Delta t < \Delta x
\label{perlastabilita}
\end{equation}
seems reasonable as constraint for stability of an explicit scheme. 

\noindent If we take as reference a $400 \times 400$ mesh, $\Delta x$ would be equal to $0.0025$.
Notwithstanding, experiments suggest the threshold $\Delta t$ should be of order $10^{-5}$ for RK-FDWENO-LF.
In the other cases, the RK-FDWENO-RK scheme, the TS-DSLWENO-LF scheme and the TS-CSLWENO-LF scheme, 
if the CFL parameter or the $\Delta t$
are adapted so as to fulfill \eqref{perlastabilita}, the simulations appear stable.

\subsection{Quality of the results}
\begin{figure}[ht]
\centering{\hspace*{-4em}  \includegraphics[scale = 1]{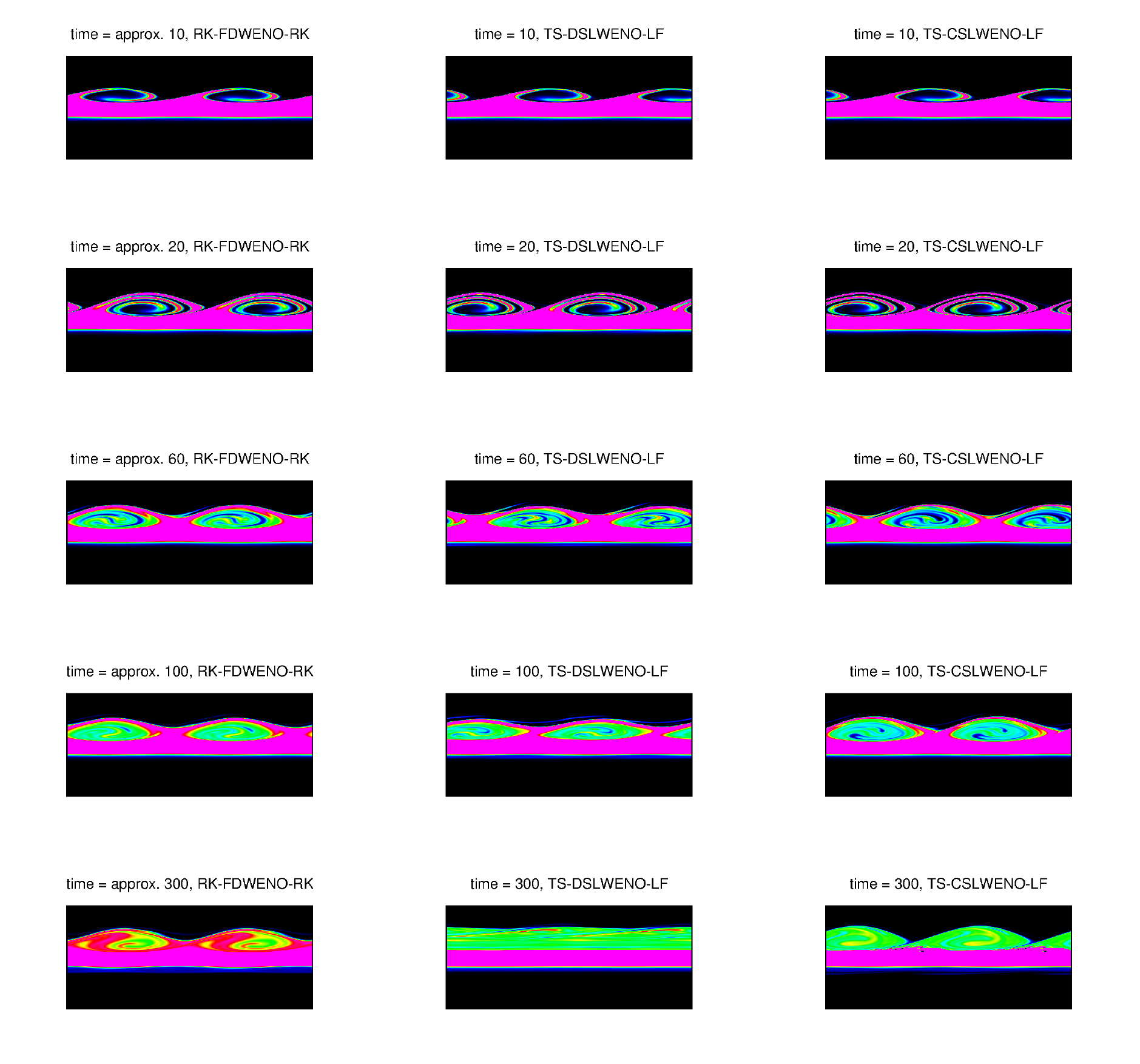} }
  \caption{
    \textbf{Comparison.} Evolution of the system up to time $\approx 300$, for a $400 \times 400$ mesh, using three different schemes. 
    Left column: the explicit conservative RK-FDWENO-RK. 
    Central column: the implicit non-conservative TS-DSLWENO-LF. 
    Right column: the implicit conservative TS-CSLWENO-LF. 
  }
  \label{evol_compared}
\end{figure}

\noindent
On Figure \ref{evol_compared} we compare at similar stages the evolution 
computed by the three most stable schemes. The dynamic of laser-plasma interaction~\cite{GBS+90,Huot2003512,BLG+08} is precisely captured.
The plasma wave, initiated by the initial fluctuations of the electron density,
exchanges energy with the electrons and with the transverse electromagnetic wave.
Vortices appear in phase space, due to the particles 
getting trapped by the plasma wave's potential well and bouncing on its separatrices. The vortices
show an oscillating behavior: they periodically inflate and deflate.
One observes the well-known ``filamentation'' phenomenon: thin
structures appear, then they are stretched thinner and folded, again
and again.  

We see that in the short term both the RK-based and the TS-based schemes behave well,
but TS-CSLWENO-LF diffuses the microscopic details more than RK-FDWENO-RK, as the long-time behavior ($t = 300$) shows.

\begin{figure}[ht]
  \centering
  \includegraphics[scale = 1]{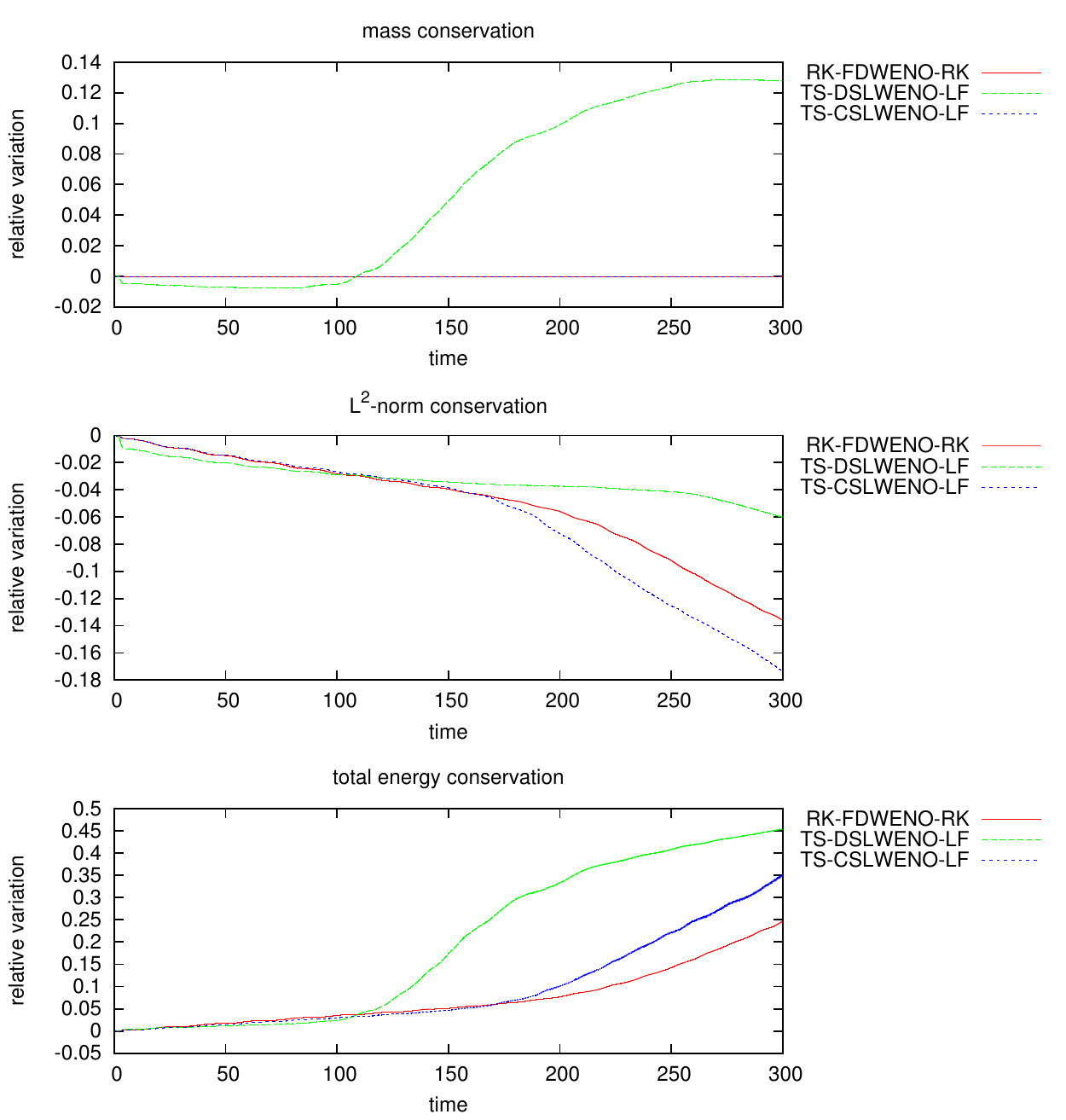}
  \caption{
    \textbf{Conservation properties.} 
    Top: the relative variation of total mass w.r.t. the initial condition.
    Center: the relative variation of the $L^2$-norm w.r.t. the initial condition.
    Bottom: the relative variation of the total energy \eqref{totalenergy} w.r.t. the initial condition.
  }
  \label{toutensemble}
\end{figure}

On Figure \ref{toutensemble} we plot the conservation properties:
the relative variation (w.r.t. time $t=0$) of the mass, of the $L^2$-norm and of the total energy
\begin{equation}
W(t)=\overbrace{\underbrace{\frac{1}{2} \int_0^1 A^2 \varrho \, \mathrm{d}x}_{\mathrm{WTK}(t) := \mathrm{kinetic}} 
+ \underbrace{\frac{\eta^{2}}{2} \int_0^1 \left[ \mathcal{E}^2 + B^2 \right] \mathrm{d}x}_{\mathrm{WTP}(t) := \mathrm{potential}}}^{\mathrm{WT}(t) := \mathrm{transversal}}
+
\overbrace{\underbrace{\int_0^1 \int_{\mathbb{R}} \left( \sqrt{1+p^2} -1 \right) f \,\mathrm{d}p \, \mathrm{d}x}_{\mathrm{WLK}(t) := \mathrm{kinetic}} + \underbrace{\frac{\eta^{2}}{2} \int_0^1 E^2 \,\mathrm{d}x}_{\mathrm{WLP}(t) := \mathrm{potential}}}^{\mathrm{WL}(t) := \mathrm{longitudinal}},
\label{totalenergy}
\end{equation}
which is shown in~\cite{CarLab06} to be conserved by the system.
We observe that around time 300 TS-DSLWENO-LF has gained about 13 \% w.r.t. the normalized mass,
which means that the plasma is strongly non-neutral, hence even the integration of the
Poisson equation becomes meaningless because the periodicity is lost.
RK-FDWENO-RK conserves better the $L^{2}$-norm, i.e. the microscopic details inside the computational domain, and the total energy.

\section{Conclusion}
\label{weno-is-good-for-you}
We have performed some preliminary tests of several WENO-based schemes to simulate the 1D quasi-relativistic Vlasov--Maxwell system, which models laser-plasma interaction.
WENO schemes, with their high accuracy and robustness to the steep gradients created by filamentation, are ideally suited to capture the dynamic of this interaction.
Indeed, our test cases have reproduced the qualitative behavior known from the literature since~\cite{GBS+90}.

To decide which schemes are more suitable for the simulation of the QRVM problem,
we tested the various combinations of Table~\ref{overall}. 
Some of them immediately appear unsatisfactory, either because they require ridiculously small
time steps, or because they are strongly non-conservative.
The two strategies which show the best behavior are RK-FDWENO-RK and TS-CSLWENO-LF, which
are both conservative; the advantage of TS-CSLWENO-LF is its implicit character and weaker
constraints on the time step, while its drawback is that 
in the long time it shows a more diffusive behavior.

{F}rom the computational point of view, WENO-based schemes have several other advantages. 
They are easily parallelizable: see for instance \cite{ManCac09} for a parallel version of RK-FDWENO.
They can be made adaptive relatively easily:
see \cite{MulVec2012} for an AMR version of TS-DSLWENO, or
\cite[and references therein]{BarMarMul2011,BaeMul2006} for an AMR version of RK-FDWENO;
the built-in computation of smoothness indicators points to the regions which have to be refined (or de-refined). 
This will be presented in a future publication.

\clearpage

\appendix

\section{Constants}
\label{Constants}
\noindent The constants involved in the dimensionless system are:
\begin{eqnarray*}
&& \eta = \frac{3}{10 \pi}, \qquad
k_{\mathrm{pump}} = 4, \qquad 
k_{\mathrm{pla}} = 2, \qquad
\omega_{0} = \sqrt{\eta^{-2} +k_{\mathrm{pump}}^2 }, \qquad
A_0 = \frac{2.5}{\omega_0},\\
&& \alpha = 0.05, \qquad
v_{\mathrm{cold}} = \sqrt{\frac{15}{511}}, \qquad
v_{\mathrm{hot}} = \sqrt{\frac{100}{511}}, \qquad
\varepsilon = \frac{\sqrt{2}}{10}, \qquad
p_{\max} = 8.
\end{eqnarray*}

\section{Constants for FBMWENO}
\label{constants_fbmweno}
\noindent If we let $x \in \left] x_{i-1/2}, x_{i+1/2} \right[$ and the interpolant is centered in the stencil, 
\begin{equation*}
\mathrm{left} = i-5/2, \qquad \mathrm{right} = i+5/2, \qquad N_{\mathrm{sub}} = 3.
\end{equation*}
The polynomials $\left\{ C_{\ell}(x) \right\}_{\ell=0}^{N_{\mathrm{sub}}-1}$ are
\begin{footnotesize}
\begin{equation*}
C_{0} =  \frac{\left( x - x_{i-5/2} \right)\left( x - x_{i-3/2} \right)}{20 \, \Delta x^{2}}, \quad
C_{1} = -\frac{\left( x - x_{i-5/2} \right)\left( x - x_{i+5/2} \right)}{10 \, \Delta x^{2}}, \quad
C_{2} =  \frac{\left( x - x_{i+5/2} \right)\left( x - x_{i+3/2} \right)}{20 \, \Delta x^{2}}.
\end{equation*}
\end{footnotesize}
The smoothness indicators are
\begin{footnotesize}
\begin{align*}
\sigma_{0} &= 
\frac{10}{3} \left( U_{i-1/2} \right)^{2} - 17 \, U_{i-1/2} \, U_{i+1/2} + 14 \, U_{i-1/2} \, U_{i+3/2} -\frac{11}{3} U_{i-1/2} \, U_{i+5/2}
+ 22 \left( U_{i+1/2} \right)^{2} \\& - \frac{111}{3} U_{i+1/2} \, U_{i+3/2} + 10 \, U_{i+1/2} \, U_{i+5/2}
+ 16 \left( U_{i+3/2} \right)^{2} - 9 \, U_{i+3/2} \, U_{i+5/2}
 +  \frac{4}{3} \left( U_{i+5/2} \right)^{2},\\[3mm]
\sigma_{1} &= 
  \frac{4}{3} \left( U_{i-3/2} \right)^{2} -7 \, U_{i-3/2} \, U_{i-1/2} + 6 \, U_{i-3/2} \, U_{i+1/2} - \frac{5}{3} U_{i-3/2} \, U_{i+3/2}
  + 10 \left( U_{i-1/2} \right)^{2} \\&- 19 \, U_{i-1/2} \, U_{i+1/2} + 6 \, U_{i-1/2} \, U_{i+3/2}
  + 10 \left( U_{i+1/2} \right)^{2} - 7 \, U_{i+1/2} \, U_{i+3/2} + \frac{4}{3} \left( U_{i+3/2} \right)^{2},\\[3mm]
\sigma_{2} &= 
\frac{10}{3} \left( U_{i+1/2} \right)^{2}
- 17 \, U_{i+1/2} \, U_{i-1/2} 
+ 14 \, U_{i+1/2} \, U_{i-3/2} 
-\frac{11}{3} U_{i+1/2} \, U_{i-5/2}
+ 22 \left( U_{i-1/2} \right)^{2}
\\&- \frac{111}{3} U_{i-1/2} \, U_{i-3/2} 
 + 10 \, U_{i-1/2} U_{i-5/2}
+ 16 \left( U_{i-3/2} \right)^{2} -9 \, U_{i-5/2} \, U_{i-3/2} +
  \frac{4}{3} \left( U_{i-5/2} \right)^{2}.
\end{align*}
\end{footnotesize}

\section*{Acknowledgments}
\noindent Francesco Vecil and Pep Mulet acknowledge financial support
from MINECO project MTM2011-22741.


\begin{thebibliography}{00}

\bibitem{GBS+90}
A. Ghizzo, P. Bertrand, M. Shoucri, T.~W. Johnston, E. Fijalkow, M.~R. Feix, 
A Vlasov code for the numerical simulation of stimulated Raman scattering, 
J. Comput. Phys.~90 (2), (1990) 431--457.

\bibitem{Huot2003512}
F.~Huot, A.~Ghizzo, P.~Bertrand, E.~Sonnendr{\"u}cker, O.~Coulaud,
       {{Instability
  of the time splitting scheme for the one-dimensional and relativistic
  Vlasov--Maxwell system}}, Journal of Computational Physics 185~(2) (2003) 512
  -- 531.

\bibitem{CarLab06}
J.~A. Carrillo, S.~Labrunie, Global solutions for the one-dimensional
  {V}lasov-{M}axwell system for laser-plasma interaction, Math. Models Methods
  Appl. Sci. 16~(1) (2006) 19--57.

\bibitem{BoCr09}
M. Bostan, N. Crouseilles
Convergence of a semi-Lagrangian scheme for the reduced Vlasov--Maxwell system for laser-plasma interaction
Numer. Math. 112 (2009), 169--195.

\bibitem{BLG+08}
N. Besse, G. Latu, A. Ghizzo, E. Sonnendr\"ucker, P.~Bertrand, 
A wavelet-MRA-based adaptive semi-Lagrangian method for the relativistic Vlasov--Maxwell system,
J. Comput. Phys.~227 (16) (2008), 7889--7916.

\bibitem{CCGMS}
M.~C\'{a}ceres, J.~Carrillo, I.~Gamba, A.~Majorana, C.-W. Shu, Deterministic
  kinetic solvers for charged particle transport in semiconductor devices,
  Cercignani, C., Gabetta, E. (eds.) Transport Phenomena and Kinetic Theory:
  Applications to Gases, Semiconductors, Photons and Biological Systems,
  Series: Modelling and Simulation in Science, Engineering and Technology,
  Birkh{\"a}user.

\bibitem{strang}
G.~Strang, On the construction and comparison of difference schemes, SIAM J.
  Numer. Anal.~(5) (1968) 506--517.

\bibitem{CarVec05}
J.~A. Carrillo, F.~Vecil, Nonoscillatory interpolation methods applied to
  {V}lasov-based models, SIAM J. Sci. Comput. 29~(3) (2007) 1179--1206
  (electronic).

\bibitem{AraBaeBelMul2011}
F.~Ar{\`a}ndiga, A.~Baeza, A.~M. Belda, P.~Mulet, {Analysis of WENO schemes for
  full and global accuracy}, SIAM Journal on Numerical Analysis 49~(2) (2011)
  893--915.

\bibitem{MR2607148}
F.~Ar{\`a}ndiga, A.~M. Belda, P.~Mulet,
       {Point-value {WENO} multiresolution applications to stable image compression}, J. Sci. Comput.
  43~(2) (2010) 158--182.

\bibitem{FilSonBer01}
F.~Filbet, E.~Sonnendr{\"u}cker, P.~Bertrand, Conservative numerical schemes
  for the {V}lasov equation, J. Comput. Phys. 172~(1) (2001) 166--187.

\bibitem{MulVec2012}
P.~Mulet, F.~Vecil,
{{A semi-Lagrangian AMR scheme for 2D transport problems in conservation form}},
  Journal of Computational Physics 237 (2013) 151--176.

\bibitem{CGMS-JCP}
J.~Carrillo, I.~Gamba, A.~Majorana, C.-W. Shu, A {WENO}-solver for the
  transients of {B}oltzmann-{P}oisson system for semiconductor devices.
  {P}erformance and comparisons with {M}onte {C}arlo methods, J. Comput.
  Phys.~(184)  498--525.

\bibitem{CGMS-JCP2}
J.~Carrillo, I.~Gamba, A.~Majorana, C.-W. Shu, 2{D} semiconductor device
  simulations by {WENO}-{B}oltzmann schemes: efficiency, boundary conditions
  and comparison to {M}onte {C}arlo methods, J. Comput. Phys.~(214) (2006)
  55--80.

\bibitem{CheKno76}
C.~Cheng, G.~Knorr, The integration of the {V}lasov equation in configuration
  space, J. Comput. Phys.~(22) (1976) 330--351.

\bibitem{JakMarOwr2000}
Z.~Jackiewicz, A.~Marthinsen, B.~Owren, {Construction of Runge--Kutta methods of
  Crouch-Grossman type of high order}, Advances in Computational Mathematics
  13~(4) (2000) 405--415.

\bibitem{MarOwr2001}
A.~Marthinsen, B.~Owren, {A Note on the Construction of Crouch--Grossman
  Methods}, BIT Numerical Mathematics 41~(1) (2001) 207--214.

\bibitem{JS}
G.-S. Jiang, C.-W. Shu, Efficient implementation of weighted {ENO} schemes, J.
  Comput. Phys. 126~(1) (1996) 202--228.

\bibitem{ManCac09}
J.~M. Mantas, M.~J. C\'aceres,
       {{Efficient
  deterministic parallel simulation of {2D} semiconductor devices based on
  {WENO-Boltzmann} schemes}}, Computer Methods in Applied Mechanics and
  Engineering 198~(5-8) (2009) 693--704.

\bibitem{BarMarMul2011}
A.~Baeza, A.~Mart\'inez-Gavara, and P.~Mulet.
\newblock Adaptation based on interpolation errors for high order mesh
  refinement methods applied to conservation laws.
\newblock {\em Applied Numerical Mathematics}, 62(4):278 -- 296, 2012.
\newblock Third Chilean Workshop on Numerical Analysis of Partial Differential
  Equations (WONAPDE 2010).

\bibitem{BaeMul2006}
A.~Baeza, P.~Mulet, 
{Adaptive mesh refinement techniques for high-order shock capturing schemes for multi-dimensional hydrodynamic simulations}, 
Int. J. Numer. Meth. Fluids 52 (2006) 455--471.

\end{thebibliography}
\end{document}